\documentclass[reqno,12pt, oneside]{amsart}

\usepackage{enumerate}
\usepackage[nobysame]{amsrefs}
\usepackage{enumitem}
\usepackage{amssymb}
\usepackage{amsmath}
\usepackage{amsthm}
\usepackage{mathtools}
\usepackage{amscd}
\usepackage{microtype}
\usepackage[right=2cm,left=2cm, top=2.5cm, bottom=2.5cm,includehead]{geometry}
\usepackage{amsfonts}
\usepackage{booktabs}
\usepackage{hyperref}

\usepackage{tgpagella}
\usepackage{eulervm}

\usepackage{chngcntr}
\counterwithin{table}{section}

\usepackage[backgroundcolor=lightgray]{todonotes}

\theoremstyle{plain}
\newtheorem{theorem}{Theorem}
\newtheorem{corollary}[theorem]{Corollary}
\newtheorem{lemma}[theorem]{Lemma}
\newtheorem{proposition}[theorem]{Proposition}

\theoremstyle{remark}
\newtheorem{remark}[theorem]{Remark}
\newtheorem{example}[theorem]{Example}

\theoremstyle{definition}

\DeclarePairedDelimiter{\abs}{\lvert}{\rvert}

\DeclarePairedDelimiter{\norm}{\lVert}{\rVert}

\DeclarePairedDelimiterX{\dset}[2]{\lbrace}{\rbrace}{#1\;\delimsize|\;#2}

\DeclareMathOperator{\var}{var}

\newcommand{\ball}{{\overline{B}}}

\DeclareMathOperator{\Lip}{Lip}
\DeclareMathOperator{\lip}{lip}
\DeclareMathOperator{\qC}{qC}

\title[Nonlinear composition operators in $bv_p$-spaces]{Nonlinear composition operators in $bv_p$-spaces:\\ continuity and compactness}

\author{Daria Bugajewska}
\address[D.~Bugajewska]{Department of Nonlinear Analysis and Applied Topology\\
Faculty of Mathematics and Computer Science\\
  Adam Mickiewicz University in Pozna\'n\\
  ul.\ Uniwersytetu Pozna\'nskiego 4\\
  61-614 Pozna\'n\\
  Poland}
\email[D.~Bugajewska]{dbw@amu.edu.pl}

\author{Piotr Kasprzak}
\address[P. Kasprzak]{Department of Nonlinear Analysis and Applied Topology\\
Faculty of Mathematics and Computer Science\\
  Adam Mickiewicz University in Pozna\'n\\
  ul.\ Uniwersytetu Pozna\'nskiego 4\\
  61-614 Pozna\'n\\
  Poland}
\email[P.~Kasprzak]{kasp@amu.edu.pl}


\date{\today}

\begin{document}
\begin{abstract}
Continuing the study initiated in our earlier article~\cite{BK}, this paper aims to characterize various continuity properties of nonlinear composition operators
acting on some sequence spaces, giving special attention to the space of sequences of bounded variation. In addition to pointwise, uniform, and locally uniform continuity, we investigate Lipschitz continuity as well as several types of H\"older continuity. Furthermore, we provide a characterization of the compactness properties of these operators.  
\end{abstract}

\maketitle

\section{Introduction}

It is well-known that the theory of nonlinear superposition operators is an essential part of the broader field of nonlinear analysis. Following~\cite{ABM}, three central and natural problems in this theory concern determining necessary and sufficient conditions on the generating function under which the superposition operator acts between given spaces, is bounded, and/or is continuous.

This paper investigates \emph{nonlinear composition operators} $C_f$ (henceforth referred to simply as \emph{composition operators}), which are superposition (or Nemytskii) operators generated by functions that do not depend on the ``time'' variable (for a precise definition, see Section~\ref{sec:composition_operators}). Our primary focus is on composition operators acting on specific sequence spaces, particularly the space $bv_p(E)$, consisting of $E$-valued sequences of $p$-bounded variation, where $p\geq 1$ and $E$ is an arbitrary normed space. This choice is motivated by the growing interest in $bv_p$-spaces and their applications (see, for example,~\cites{BA, BAM, K, LF} and references therein). Additionally, we aim to extend the study of Dedagich and Zabre\u{\i}ko (see~\cite{DZ}) on the topological properties of superposition operators in classical sequence spaces -- especially in $l^p$-spaces -- to a more general setting. In this context, one can also mention the paper~\cite{KFA}, which is the first to discuss the properties of superposition operators in the case of $bv_p$-spaces. However, some concerns have been raised regarding its results (for a discussion of these issues, see \cite{BK}).  

In the paper~\cite{BK}, we provided a complete description of the acting conditions for composition operators in several classical sequence spaces, including $c_0(E)$, $c(E)$, $l^p(E)$, as well as $bv_p$-spaces. In addition, we characterized their boundedness and local boundedness. In the present article, we continue our study by investigating the topological properties of these operators. We focus on describing different types of continuity, thus addressing the third problem outlined above. In particular, we formulate conditions on the generator $f$ under which the operator $C_f$ is pointwise continuous, (locally) uniformly continuous, and H\"older continuous.  

Another interesting problem concerns the compactness of a composition operator. It is well-known that such an operator can only be compact in various function spaces -- such as spaces of continuous or integrable functions, or functions of bounded variation -- if its generator is constant. As demonstrated in this article, a similar result holds for composition operators acting between different sequence spaces.

The paper is organized as follows. In Section 2, we introduce the necessary notations and recall basic results concerning classical sequence spaces. We also define several classes of H\"older continuous mappings, highlighting their key properties and relationships among them. Given the length of this paper and the accessibility of our earlier work~\cite{BK}, we have opted for a concise treatment of the preliminaries, avoiding repetition of conventions, motivations, and previously established results concerning the H\"older classes. In Section~3, we study the continuity properties of composition operators acting on various sequence spaces. In particular, we examine pointwise continuity, uniform continuity, and local uniform continuity, that is, uniform continuity on bounded sets. In Section 4, we provide both necessary and sufficient conditions on the generating function under which composition operators exhibit Lipschitz and H\"older continuity, with a clear distinction between global and local settings. Section 5 is devoted to the study of the compactness of composition operators. For the reader’s convenience, each section concludes with one or more tables that summarize the main results in a clear and accessible format. In compiling these tables, we have also incorporated relevant results on acting conditions and boundedness of composition operators in sequence spaces established in~\cite{BK}. Finally, in Section~6 we compare our findings with corresponding theorems from the parallel theory of composition operators on spaces of functions of Wiener bounded $p$-variation, and provide references to those results.

\section{Preliminaries}

In this short section, we introduce notation and recall some basic facts related to classical sequence spaces and composition operators. Our exposition closely follows that of~\cite{BK}, and we focus on the most relevant aspects here. For a more detailed treatment, we refer the reader to~\cite{BK}*{Section~2}. 

Throughout the paper, let $(E,\norm{\cdot})$, or simply $E$, denote a normed space over either the real or complex field. By $B_E(u,r)$ and $\ball_E(u, r)$ we denote 
the open and closed balls in $E$, centered at $u \in E$ with radius $r > 0$, respectively. 

\subsection{Sequence spaces} Elements of sequence spaces will be denoted using the functional notation; that is, $x(n)$ will represent the $n$-th term of the sequence $x$. In contrast, we will use the classical notation $(u_n)_{n \in \mathbb N}$ when referring to sequences of other objects, possibly elements of sequence spaces. 

In this paper, we consider classical sequence spaces, including $l^p(E)$ for $p \in [1,+\infty)$, which consists of all $E$-valued sequences that are absolutely summable with $p$-th power. This space is equipped with the norm
\[
 \norm{x}_{l^p} := \Biggl(\sum_{n=1}^\infty \norm{x(n)}^p \Biggr)^{\frac{1}{p}}.
\]
We also consider $l^\infty(E)$, the normed space of all bounded $E$-valued sequences, endowed with the supremum norm
\[
 \norm{x}_\infty := \sup_{n \in \mathbb N}\norm{x(n)}.
\]
Additionally, we consider its normed subspaces: $c(E)$, consisting of all convergent sequences, and $c_0(E)$, the space of null sequences. The primary focus of our study, however, is the space $bv_p(E)$ of all sequences $x \colon \mathbb N \to E$ such that
\[
 \sum_{n=1}^\infty \norm{x(n+1)-x(n)}^p<+\infty.
\]
The standard norm on $bv_p(E)$ is given by 
\[
 \norm{x}_{bv_p}:=\norm{x(1)}+ \Biggl(\sum_{n=1}^\infty \norm{x(n+1)-x(n)}^p\Biggr)^{\frac{1}{p}}. 
\]
If $E$ is a Banach space, then all the sequence spaces considered here are complete. 

Let us also recall inclusions holding between the above sequence spaces. For $1 \leq p < q <+\infty$ we have
\[
 l^p(E) \subset l^q(E) \subset c_0(E) \subset c(E) \subset l^\infty(E) \quad \text{and} \quad l^p(E) \subset bv_p(E) \subset bv_q(E). 
\]
Additionally, $bv_1(E) \subset l^{\infty}(E)$. However, for $p>1$ neither $bv_p(E)$ is contained in $l^{\infty}(E)$, nor $l^{\infty}(E)$ is contained in $bv_p(E)$. In particular,for $p>1$ we have $bv_p(E) \not\subseteq c(E)$. Finally, if $E$ is a complete normed space, then $bv_1(E)\subset c(E)$ and $\norm{x}_\infty \leq \norm{x}_{bv_1}$ for any $x \in bv_1(E)$. For more details or proofs regarding these inclusions for $bv_p$-spaces in the case $E:=\mathbb R$, we refer the reader to~\cite{BA}. 

\subsection{Composition operators and their generators}
\label{sec:composition_operators}

Let  $f \colon E \to E$ be a map. A \emph{composition operator} (or, more precisely, a \emph{nonlinear composition operator}) $C_f$, associated with the function $f$, assigned to a given sequence $x$ (which usually belongs to some specific set or space) the sequence $C_f(x)$ whose $n$-th term is defined by $C_f(x)(n)=f(x(n))$. The function $f$ is referred to as the \emph{generator} of $C_f$. For further details on composition and superposition operators, we refer the reader to the celebrated monograph by Appell and Zabrejko~\cite{AZ}.

Natural classes of maps generating composition operators in $bv_p$-type spaces consists of H\"older continuous mappings. Since there is no consensus in the literature regarding the naming of the various classes of such functions, we provided a comprehensive overview of the necessary definitions in~\cite{BK}, along with a discussion of their properties and interrelationships. For the sake of clarity and the reader's convenience, we briefly recall these definitions below.

Let  $\alpha\in (0,1]$ and $p,q \in [1,+\infty)$. A map $f\colon E \to E$ is called
\begin{enumerate}[label=\textup{(\alph*)}]
 \item \emph{H\"older continuous with exponent $\alpha$}, if there exists $L\geq 0$ such that $\norm{f(u)-f(w)}\leq L\norm{u-w}^\alpha$ for any $u,w \in E$; the set of all such maps will be denoted by $\Lip^\alpha(E)$,

 \item \emph{locally H\"older continuous in the stronger sense with exponent $\alpha$}, if there exist $\delta>0$ and $L\geq 0$ such that $\norm{f(u)-f(w)}\leq L\norm{u-w}^\alpha$ for all $u,w \in E$ with $\norm{u-w}\leq \delta$; the set of all such maps will be denoted by $\Lip^\alpha_{\text{strg}}(E)$,
  
 \item\label{local_hoder} \emph{H\"older continuous on bounded sets with exponent $\alpha$}, if for every $r>0$ there exists $L_r \geq 0$ such that $\norm{f(u)-f(w)}\leq L_r\norm{u-w}^\alpha$ for all $u,w \in \ball_E(0,r)$; the set of all such maps will be denoted by $\Lip^\alpha_{\text{bnd}}(E)$,	

 \item \emph{H\"older continuous on compact sets with exponent $\alpha$}, if for every compact subset $K$ of $E$ there exists $L_K \geq 0$ such that $\norm{f(u)-f(w)}\leq L_K\norm{u-w}^\alpha$ for all $u,w \in K$; the set of all such maps will be denoted by $\Lip^\alpha_{\text{comp}}(E)$,

 \item\label{it:little_H} \emph{locally little H\"older continuous with exponent $\alpha$}, if for every $r>0$ 
\[
 \lim_{\delta\to 0^+} \sup\dset[\Bigg]{\frac{\norm{f(u)-f(w)}}{\norm{u-w}^\alpha}}{\text{$u,w \in \ball_E(0,r)$ with $0<\norm{u-w} \leq \delta$}} = 0,
\]
or equivalently, if for every $r>0$ and every $\varepsilon>0$ there exists $\delta>0$ such that $\norm{f(u)-f(w)}\leq \varepsilon \norm{u-w}^\alpha$ for all $u,w \in \ball_E(0,r)$ with $\norm{u-w} \leq \delta$; the set of all such maps will be denoted by $\lip^\alpha_{\text{loc}}(E)$.
\end{enumerate}

\noindent When the exponent $\alpha$ is clear from the context and there is no risk of confusion, we will omit the phrase ``with exponent'' and simply write, for example, ``$f$ is H\"older continuous''  instead of ``$f$ is H\"older continuous with exponent $\alpha$.'' Additionally, maps that are H\"older continuous (in any sense) with exponent $\alpha=1$ will be referred to as Lipschitz continuous.

When discussing the local uniform continuity of composition operators, we will encounter yet another condition related to the H\"older continuity of the generators. We have chosen to present it separately, as it is well-defined only in the case $E:=\mathbb R$. Let $\alpha \in (0,1]$. A map $f \colon \mathbb R \to \mathbb R$ is called  
\begin{enumerate}[label=\textup{(\alph*)}, resume]
 \item\label{it:g} \emph{quasi continuously differentiable with exponent $\alpha$}, if for every $r>0$ and $\varepsilon>0$ there exists $\delta>0$ such that for all $a,b,c,d \in [-r,r]$ with $a\neq b$, $c\neq d$ and $\abs{a-c}\leq \delta$, $\abs{b-d}\leq \delta$ the following inequality holds
\begin{equation}\label{eq:quasiC1}
 \abs[\Bigg]{\dfrac{(f(a)-f(b))(a-b)}{\abs{a-b}^{1+\alpha}} - \dfrac{(f(c)-f(d))(c-d)}{\abs{c-d}^{1+\alpha}} } \leq \varepsilon;
\end{equation} 
the set of all such functions will be denoted by $\qC^{1,\alpha}(\mathbb R)$.
\end{enumerate}

\begin{remark}
The form of condition~\ref{it:g} certainly requires clarification. A seemingly more natural and simpler alternative would be to express~\eqref{eq:quasiC1} as
\begin{equation}\label{eq:quasiC2}
 \abs[\Bigg]{\dfrac{f(a)-f(b)}{\abs{a-b}^{\alpha}} - \dfrac{f(c)-f(d)}{\abs{c-d}^{\alpha}} } \leq \varepsilon.
\end{equation} 
Then, however, for any two points $a,b \in [-r,r]$ with $0<\abs{a-b}\leq \delta$, substituting both $(a,b)$ and $(b,a)$ into the left-hand side of~\eqref{eq:quasiC2} would yield 
\[
 \varepsilon \geq \abs[\Bigg]{\dfrac{f(a)-f(b)}{\abs{a-b}^{\alpha}} - \dfrac{f(b)-f(a)}{\abs{b-a}^{\alpha}} } = \frac{2\abs{f(a)-f(b)}}{\abs{a-b}^\alpha}.
\]  
This, in turn, would imply that $f$ satisfies the local little H\"older condition, which is a stronger requirement than what we have intended to impose (see Proposition~\ref{prop:little_vs_qC1} below). In contrast, substituting $(a,b)$ and $(b,a)$ into the left-hand side of~\eqref{eq:quasiC1} results in zero.
\end{remark}

The discussion above suggests a potential connection between conditions~\ref{it:little_H} and~\ref{it:g}. To explore this relationship further, we introduce an auxiliary map $\Phi_\alpha$, where $\alpha \in (0,1]$. Given a function $f \colon \mathbb R \to \mathbb R$, we define $\Phi_\alpha(f)$ by the formula
\[
 \Phi_\alpha(f)(u,w):=\dfrac{(f(u)-f(w))(u-w)}{\abs{u-w}^{1+\alpha}}\ \ \text{for $u\neq w$.}
\]
Clearly, a function $f$ is quasi continuously differentiable with exponent $\alpha$ if and only if $\Phi_\alpha(f)$ is uniformly continuous on $\Delta_r:=([-r,r]\times [-r,r])\setminus \dset{(x,x)}{x \in \mathbb R}$ for every $r>0$. Moreover, note that if $f \colon \mathbb R \to \mathbb R$ is continuous, then the function $\Phi_\alpha(f)$ is also continuous.

With this correspondence between $f$ and $\Phi_\alpha(f)$ at hand, we can readily establish that $\lip^\alpha_{\text{loc}}(\mathbb R) \subseteq \qC^{1,\alpha}(\mathbb R)$.

\begin{proposition}\label{prop:little_vs_qC1}
Let $\alpha \in (0,1]$. If a function $f \colon \mathbb R\to \mathbb R$ is locally little H\"older continuous with exponent~$\alpha$, then it is also quasi continuously differentiable with the same exponent.
\end{proposition}

\begin{proof}
Fix $r>0$ and $\varepsilon>0$. Since $f \in \lip^\alpha_{\text{loc}}(\mathbb R)$, there exists $\delta>0$ such that for all $u,w \in [-r,r]$ satisfying $\abs{u-w} \leq \delta$ we have $\abs{f(u)-f(w)}\leq \varepsilon\abs{u-w}^\alpha$.  Consequently, for distinct $u,w \in [-r,r]$ with $\abs{u-w}\leq \delta$, it follows that
\[
 \abs{\Phi_\alpha(f)(u,w)} \leq \frac{\abs{f(u)-f(w)}\cdot \abs{u-w}}{\abs{u-w}^{1+\alpha}}\leq \varepsilon.
\]
Putting $\Phi_\alpha(f)(u,u):=0$ for $u \in \mathbb R$, we extend the function $\Phi_\alpha(f)$ over the entire square $[-r,r]\times [-r,r]$ to a uniformly continuous map. This, in turn, implies that $f$ belongs to $\qC^{1,\alpha}(\mathbb R)$. 
\end{proof}

The reverse implication does not hold when $\alpha=1$, as $\lip^1_{\text{loc}}(\mathbb R)$ consists only of constant functions, whereas $\qC^{1,1}(\mathbb R)$ includes, for example, the identity function. (In fact, as we will see shortly, $\qC^{1,1}(\mathbb R)$ coincides with $C^1(\mathbb R)$, that is, the space  of continuously differentiable functions.) For $\alpha \in (0,1)$ it remains an open question whether the classes $\lip^\alpha_{\text{loc}}(\mathbb R)$ and $\qC^{1,\alpha}(\mathbb R)$ are distinct.

Next, we explore the relationship between $\qC^{1,\alpha}(\mathbb R)$ and $\Lip^\alpha_{\text{bnd}}(\mathbb R)$.

\begin{proposition}\label{prop:qc1_vs_lip_loc}
Let $\alpha \in (0,1]$. If a function $f \colon \mathbb R \to \mathbb R$ is quasi continuously differentiable with exponent $\alpha$, then it is also H\"older continuous on bounded sets with the same exponent. 
\end{proposition}

\begin{proof}
It suffices to show that the function $\Phi_\alpha(f)$ is bounded on each set $\Delta_\rho$ for $\rho>0$. Suppose, for the sake of contradiction, that there exists a radius $R>0$ and a sequence $(u_n,w_n)_{n \in \mathbb N}$ in $\Delta_R$ such that $\abs{\Phi_\alpha(f)(u_n,w_n)} \to +\infty$ as $n \to +\infty$. Since the set $\Delta_R$ is relatively compact in $[-R,R]\times [-R,R]$, the sequence $(u_n,w_n)_{n \in \mathbb N}$ admits a Cauchy subsequence $(u_{n_k},w_{n_k})_{k \in \mathbb N}$. In particular, there exists an index $m$ such that $\abs{u_{n_m} - u_{n_l}} \leq \delta$ and $\abs{w_{n_m} - w_{n_l}} \leq \delta$ for all $l \geq m$; here, $\delta>0$ is the constant from condition~\ref{it:g} corresponding to $\varepsilon:=1$ and $r:=R$. Then, for all $l \geq m$ we have
\begin{align*}
 \abs[\big]{\Phi_\alpha(f)(u_{n_l},w_{n_l})} & \leq \abs[\big]{\Phi_\alpha(f)(u_{n_l},w_{n_l})-{\Phi_\alpha(f)(u_{n_m},w_{n_m})}} + \abs[\big]{\Phi_\alpha(f)(u_{n_m},w_{n_m})}\\
 & \leq 1 + \abs[\big]{\Phi_\alpha(f)(u_{n_m},w_{n_m})}. 
\end{align*}
This implies that $\sup_{l \geq m}\abs[\big]{\Phi_\alpha(f)(u_{n_l},w_{n_l})}<+\infty$, contradicting our initial assumption that $\abs{\Phi_\alpha(f)(u_n,w_n)} \to +\infty$ as $n \to +\infty$. Therefore, $\Phi_\alpha(f)$ is bounded on each set $\Delta_\rho$, and consequently, $f$ is H\"older continuous on bounded sets with exponent $\alpha$.
\end{proof}

The fact that, in general, the classes $\qC^{1,\alpha}(\mathbb R)$ and $\Lip^\alpha_{\text{bnd}}(\mathbb R)$ do not coincide can be illustrated by the following simple example.

\begin{example}
Let $\alpha \in (0,1]$, and consider the power function $f \colon \mathbb R \to \mathbb R$ defined by $f(u)=\abs{u}^\alpha$. It is well-known that $f$ satisfies the H\"older condition (on bounded sets) with exponent $\alpha$. However, for any $\delta>0$ and $0<u\leq \frac{1}{2}\delta$, we have
\begin{align*}
  \abs[\Bigg]{\dfrac{(f(u)-f(0))(u-0)}{\abs{u-0}^{1+\alpha}} - \dfrac{(f(-u)-f(0))(-u-0)}{\abs{-u-0}^{1+\alpha}} }  =2.
\end{align*}
This implies that $f$ does not belong to $\qC^{1,\alpha}(\mathbb R)$.
\end{example}

The term ``quasi continuously differentiable'' suggests that condition~\ref{it:g} should have some connection to continuous differentiability. Indeed, this is the case. It turns out that for $\alpha=1$ those two notions coincide. 

\begin{proposition}\label{prop:qc1_vs_c1}
A function $f \colon \mathbb R \to \mathbb R$ is continuously differentiable if and only if it is quasi continuously differentiable with exponent $\alpha=1$.  
\end{proposition}   

\begin{proof}
The proof of the right implication proceeds similarly to the proof of Proposition~\ref{prop:little_vs_qC1}. Let $r>0$ be fixed. We define the function $\Psi(f) \colon [-r,r]\times [-r,r] \to \mathbb R$ as
\[
 \Psi(f)(u,w):=\begin{cases} 
           \dfrac{f(u)-f(w)}{u-w}, & \text{if $u\neq w$,}\\[2mm]
           f'(u), & \text{if $u=w$.} 
			   \end{cases}
\]  
Since $f$ is continuously differentiable on $(-2r,2r)\times (-2r,2r)$, the function $\Psi(f)$ is uniformly continuous on its domain. Furthermore, for all distinct $u,w \in [-r,r]$ we have
\[
 \Phi_1(f)(u,w)= \frac{(f(u)-f(w)) (u-w)}{(u-w)^2}=\frac{f(u)-f(w)}{u-w}=\Psi(f)(u,w).
\]
This shows that $\Phi_1(f)$ is uniformly continuous on $\Delta_r$. As a result, $f$ is quasi continuously differentiable with exponent $\alpha=1$.

The proof of the opposite implication relies on a result concerning absolutely continuous maps. It states that if a real-valued function $g$ is absolutely continuous on an interval $[a,b]$ and $D_g \subseteq [a,b]$ denotes the set of points where $g$ is differentiable, then $g$ is continuously differentiable on $[a,b]$ if and only if the restriction of $g'$ to $D_g$ is uniformly continuous (see~\cite{reinwand}*{Lemma~1.1.27~(a)}).

Assume that $f$ is quasi continuously differentiable with exponent $\alpha=1$, and let $r>0$. By Proposition~\ref{prop:qc1_vs_lip_loc} the function $f$ is Lipschitz continuous (and thus absolutely continuous) on the interval $[-r,r]$. For a fixed $\varepsilon>0$ let $\delta>0$ be the constant from condition~\ref{it:g} corresponding to $2r$. Consider any points $u,w \in D_f$ such that $\abs{u-w}\leq \delta$; here, $D_f$ denotes the set of points in the interval $[-r,r]$ at which $f$ is differentiable. Applying~\eqref{eq:quasiC1} with $a:=u+h$, $b:=u$, $c:=w+h$, $d:=w$, where $0<h<r$, we obtain
\begin{align*}
    \abs[\Bigg]{\dfrac{f(u+h)-f(u)}{h} - \dfrac{f(w+h)-f(w)}{h}} \leq \varepsilon. 
\end{align*}
Taking the limit as $h \to 0^+$ in the above inequality, we conclude that $\abs{f'(u)-f'(w)}\leq \varepsilon$. This shows that the restriction of $f'$ to the set $D_f$ is uniformly continuous. Therefore, $f$ is continuously differentiable on $[-r,r]$.
\end{proof}

For $\alpha \in (0,1)$ the classes $C^1(\mathbb R)$ and $\qC^{1,\alpha}(\mathbb R)$ are distinct, but they are still closely related.

\begin{proposition}
If a function $f \colon \mathbb R \to \mathbb R$ is continuously differentiable, then it is also quasi continuously differentiable with any exponent $\alpha \in (0,1)$.
\end{proposition}

\begin{proof}
If $f$ is continuously differentiable, then it is obviously Lipschitz continuous on bounded sets. This immediately implies that $f \in \lip^{\alpha}_{\text{loc}}(\mathbb R)$. To complete the proof, it suffices to apply Proposition~\ref{prop:little_vs_qC1}.
\end{proof}

The fact that the inclusion $C^1(\mathbb R) \subseteq \qC^{1,\alpha}(\mathbb R)$ is strict for $\alpha \in (0,1)$ is demonstrated by the following example.

\begin{example}
Fix $\alpha \in (0,1)$ and choose $\beta$ such that $0<\alpha<\beta<1$. Define the function $f \colon \mathbb R \to \mathbb R$ by $f(u)=\abs{u}^\beta$. Then, $f \in \Lip^\beta_{\text{bnd}}(\mathbb R)$. This, in turn, implies that $f \in \lip^\alpha_{\text{loc}}(\mathbb R)$. By Proposition~\ref{prop:little_vs_qC1} we conclude that $f \in \qC^{1,\alpha}(\mathbb R)$. However, the function $f$ is not continuously differentiable on any open interval containing zero.
\end{example}

This concludes our discussion of the interrelationships between various classes of H\"older continuous maps. As these maps are not the primary focus of our study, we have decided not delve deeply into the details. For more results in this direction, we refer the reader to Section~2.5 of our earlier paper~\cite{BK}. 

\section{Continuity}

Our first goal is to examine the continuity properties of composition operators on certain sequence spaces. The concepts of pointwise and uniform continuity are well understood and require no introduction. However, it is useful to recall that a (nonlinear) operator between two normed spaces $X$ and $Y$ is \emph{locally uniformly continuous} if it is uniformly continuous on every (closed or open) ball of $X$, or equivalently, on any bounded subset of $X$.

Before we continue, let us first prove the following simple but general result.

\begin{proposition}\label{prop:cont_Cf_implies_f}
Let $p,q \in [1,+\infty)$. Moreover, let $X \in \{l^p(E),c_0(E),c(E),l^\infty(E),bv_p(E)\}$ and $Y \in \{l^q(E),c_0(E),c(E),l^\infty(E),bv_q(E)\}$. Assume also that the composition operator $C_f$ maps $X$ into $Y$. If $C_f$ is continuous \textup(respectively, locally uniformly continuous or uniformly continuous\textup), then the functin $f$ possesses the same continuity property.
\end{proposition}

\begin{proof}
(a) \emph{Pointwise continuity}. Assume that $C_f \colon X \to Y$ is continuous. Let $u \in E$ and consider any sequence $(u_n)_{n \in \mathbb N}$ in $E$ that converges to $u$. Define $x:=(u,0,0,0,\ldots)$ and $x_n:=(u_n,0,0,0,\ldots)$ for $n \in \mathbb N$. As $\norm{x-x_n}_X \leq 2\norm{u-u_n}$ for all $n \in \mathbb N$, it follows that $x_n \to x$ in $X$. Therefore, $C_f(x_n) \to C_f(x)$ in $Y$. Since $\norm{f(u_n)-f(u)}\leq \norm{C_f(x_n)-C_f(x)}_Y$, this implies that $f(u_n) \to f(u)$ in $E$. Thus, $f$ is continuous.

(b) \emph{Local uniform continuity}. Assume that $C_f \colon X \to Y$ is locally uniformly continuous. Fix $r>0$ and $\varepsilon>0$. By assumption, there exists $\delta>0$ such that $\norm{C_f(\xi)-C_f(\eta)}_Y \leq \varepsilon$ for all $\xi,\eta \in \ball_X(0,r)$ with $\norm{\xi-\eta}_X\leq \delta$. Now, let us take any $u,w \in \ball_E(0,r)$ such that $\norm{u-w}\leq \delta$. Define $x:=(u,0,0,0,\ldots)$, $y:=(w,0,0,0,\ldots)$ if $X\in \{l^p(E),c_0(E),c(E),l^\infty(E)\}$ and $x:=(u,u,u,u,\ldots)$, $y:=(w,w,w,w,\ldots)$ if $X=bv_p(E)$. In both cases, we have $x,y \in \ball_X(0,r)$ and $\norm{x-y}_X=\norm{u-w}\leq \delta$. Consequently, $\norm{f(u)-f(w)}\leq \norm{C_f(x)-C_f(y)}_Y \leq \varepsilon$. This proves that $f$ is locally uniformly continuous.

(c) \emph{Uniform continuity}. The proof is almost identical to the one presented in part~(b), so we will omit it here.
\end{proof}

\subsection{Pointwise continuity}
\label{sec:pointwise_continuity}

In the result above, we have seen that the continuity properties of the composition operator are reflected in the properties of its generator. Interestingly, in most cases, the converse is also true.

\begin{theorem}\label{thm:c_0_bvp_continuous}
Let $p \in [1,+\infty)$. Moreover, assume that the composition operator $C_f$ maps $c_0(E)$ into $bv_p(E)$. Then, $C_f$ is continuous if and only if $f$ is.
\end{theorem}

\begin{proof}
By Proposition~\ref{prop:cont_Cf_implies_f} it is enough to prove the sufficiency part. Assume that $f$ is continuous. Since $C_f$ maps $c_0(E)$ into $bv_p(E)$, it follows that $f$ is locally constant at $0$, meaning there exists $\delta>0$ such that the restriction $f|_{\ball_E(0,\delta)}$ is a constant function (see \cite{BK}*{Theorem 15}). Now, fix $\varepsilon>0$ and let $x \in c_0(E)$. Then, there is an index $N \in \mathbb N$ such that $\norm{x(n)}\leq \frac{1}{2}\delta$ for $n >N$. Furthermore, the continuity of $f$ ensures the existence of a number $\eta \in (0,\frac{1}{2}\delta]$ such that for every $n\in \{1,\ldots,N\}$ and all $w \in \ball_E(x(n),\eta)$ we have $\norm{f(x(n))-f(w)}\leq 2^{-(\frac{1}{p}+2)}N^{-\frac{1}{p}}\varepsilon$. Now, let $y \in c_0(E)$ satisfy $\norm{x-y}_\infty \leq \eta$. Then, $\norm{y(n)}\leq \delta$ for $n>N$. Consequently,
\begin{align*}
& \sum_{n=1}^\infty \norm[\Big]{\bigl[f(x(n+1))-f(y(n+1))\bigr] - \bigl[f(x(n))-f(y(n))\bigr]}^p \\
&\qquad \leq 2^p \cdot \sum_{n=1}^\infty \norm[\big]{f(x(n+1))-f(y(n+1))}^p + 2^p \cdot \sum_{n=1}^\infty \norm[\big]{f(x(n))-f(y(n))}^p\\
&\qquad \leq 2^{p+1}  \cdot \sum_{n=1}^N \norm[\big]{f(x(n))-f(y(n))}^p\\
&\qquad \leq 2^{p+1} \cdot N \cdot 2^{-(2p+1)}N^{-1}\varepsilon^p = 2^{-p}\varepsilon^p.
\end{align*}
Since $\norm{f(x(1))-f(y(1))}\leq \frac{1}{4}\varepsilon$, it follows that $\norm{C_f(x)-C_f(y)}_{bv_p}\leq \varepsilon$. Thus, $C_f$ is continuous at $x$.
\end{proof}

A similar result holds also for composition operators mapping  $l^p(E)$ into $bv_q(E)$, though the proof is slightly more technical yet follows the same overall approach.

\begin{theorem}
Let $p,q \in [1,+\infty)$. Moreover, assume that the composition operator $C_f$ maps $l^p(E)$ into $bv_q(E)$. Then, $C_f$ is continuous if and only if $f$ is .
\end{theorem}

\begin{proof}
By Proposition~\ref{prop:cont_Cf_implies_f} we may assume that $f$ is continuous. Since $C_f$ maps $l^p(E)$ into $bv_q(E)$, there exist $\delta>0$ and $L\geq 0$ such that $\norm{f(u)-f(w)}\leq L\norm{u}^{\frac{p}{q}}+L\norm{w}^{\frac{p}{q}}$ for all $u,w \in \ball_E(0,\delta)$ (see \cite{BK}*{Theorem 17}). Fix $\varepsilon>0$ and $x \in l^p(E)$. Choose a positive integer $N \geq 2$ such that $\norm{x(n)}\leq \frac{1}{2}\delta$ for $n \geq N$ and
\[
\sum_{n=N}^\infty \norm{x(n)}^p \leq \frac{\varepsilon^q}{2^{3q+p+3} (L+1)^q}.
\]
For simplicity, let $\gamma$ denote the right-hand side of the above estimate. For $n=1,\ldots,N$ let $0<\eta \leq \min\{\frac{1}{2}\delta,\gamma^{\frac{1}{p}}\}$ be such that
\[
\norm{f(x(n))-f(w)} \leq \frac{\varepsilon}{2^{2+\frac{2}{q}}N^{\frac{1}{q}}},
\]
whenever $\norm{w - x(n)}\leq \eta$. Now, take $y \in l^p(E)$ with $\norm{x-y}_{l^p}\leq \eta$. Then,
\[
  \norm{f(x(1))-f(y(1))} \leq \frac{\varepsilon}{2^{2+\frac{2}{q}}N^{\frac{1}{q}}} \leq \frac{1}{4}\varepsilon.
\]
Additionally, we have the estimate
\begin{align*}
& \sum_{n=1}^{N-1} \norm[\Big]{ \bigl[f(x(n+1))-f(y(n+1))\bigr]-\bigl[f(x(n))-f(y(n))\bigr]}^q\\
&\qquad \leq 2^{q+1} \sum_{n=1}^N \norm[\big]{ f(x(n))-f(y(n))}^q \leq 2^{q+1}N \cdot \frac{\varepsilon^q}{2^{2q+2}N} \leq \frac{1}{2^{q+1}}\varepsilon^q. 
\end{align*}
We now turn to the sum
\[
  \sum_{n=N}^{\infty} \norm[\Big]{ \bigl[f(x(n+1))-f(y(n+1))\bigr]-\bigl[f(x(n))-f(y(n))\bigr]}^q.
\]
First, let us note that $\norm{y(n)}\leq \delta$ for $n \geq N$. Also,
\[
\Biggl( \sum_{n=N}^\infty \norm{y(n)}^p \Biggr)^{\frac{1}{p}} \leq \Biggl( \sum_{n=N}^\infty \norm{y(n)-x(n)}^p \Biggr)^{\frac{1}{p}} + \Biggl( \sum_{n=N}^\infty \norm{x(n)}^p \Biggr)^{\frac{1}{p}} \leq \eta + \gamma^{\frac{1}{p}}\leq 2\gamma^{\frac{1}{p}}.
\] 
Therefore, 
\begin{align*}
& \sum_{n=N}^{\infty} \norm[\Big]{ \bigl[f(x(n+1))-f(y(n+1))\bigr]-\bigl[f(x(n))-f(y(n))\bigr]}^q\\
&\qquad \leq 2^{q+1} \sum_{n=N}^{\infty} \norm[\big]{ f(x(n))-f(y(n))}^q \leq 2^{q+1}L^q \sum_{n=N}^{\infty}\bigl(\norm{x(n)}^{\frac{p}{q}} + \norm{y(n)}^{\frac{p}{q}}\bigr)^q\\
&\qquad \leq 2^{2q+1}L^q \sum_{n=N}^{\infty} \bigl(\norm{x(n)}^p + \norm{y(n)}^p\bigr) \leq  2^{2q+1}L^q (\gamma+2^p\gamma)\leq 2^{2q+p+2}L^q\gamma\\
& \qquad =  2^{2q+p+2}L^q \cdot \frac{\varepsilon^q}{2^{3q+p+3} (L+1)^q} \leq \frac{1}{2^{q+1}}\varepsilon^q.
\end{align*}
Summing all the estimates, we conclude that
\[
 \norm{C_f(x)-C_f(y)}_{bv_q} \leq \frac{1}{4}\varepsilon + \Biggl(\frac{1}{2^{q+1}}\varepsilon^q+\frac{1}{2^{q+1}}\varepsilon^q \Biggr)^{\frac{1}{q}} \leq \frac{1}{4}\varepsilon + \frac{1}{2}\varepsilon\leq \varepsilon.
\]
Thus, the composition operator $C_f$ is continuous at $x \in l^p(E)$.
\end{proof}

Now, we shift our focus to composition operators acting between $bv_p$-spaces. In this case, we do not need to assume additionally that the generators are continuous. This is because their continuity follows directly from the general requirement that $C_f$ map one $bv_p$-space into another.

\begin{theorem}\label{thm:cont_bv_p_bv_q}
Let $p,q \in [1,+\infty)$, and let $E$ be a Banach space when $p=1$, or a normed space when $p>1$. Every composition operator $C_f$ that maps $bv_p(E)$ into $bv_q(E)$ is continuous.
\end{theorem}

\begin{proof}
The proof proceeds in three steps.
\begin{enumerate}[label=\textup{(\alph*)}, wide]

\item If $1\leq q < p < +\infty$, there is nothing to prove. In this case, the composition operator $C_f$ is a constant map, and therefore continuous (cf.~\cite{BK}*{Corollary~19}).

\item Let $p=1$ and $1 \leq q < +\infty$. Since $C_f$ maps $bv_1(E)$ into $bv_q(E)$, for each compact subset $K$ of $E$ there exists a constant $L_K \geq 0$ such that $\norm{f(u)-f(w)}\leq L_K\norm{u-w}^{\frac{1}{q}}$ for all $u,w \in K$ (see~\cite{BK}*{Theorem~20}). We aim to show that this condition implies that for every $u \in E$ there are constants $r>0$ and $L_r\geq 0$ such that $\norm{f(w)-f(v)}\leq L_r\norm{w-v}^{\frac{1}{q}}$ for all $w,v \in \ball_E(u,r)$. Suppose, for the sake of contradiction, that this is not true. That is, there exist $u \in E$ and $w_n,v_n \in \ball_E(u,\frac{1}{n})$ such that $\norm{f(v_n)-f(w_n)}>n\norm{v_n-w_n}^{\frac{1}{q}}$ for all $n \in \mathbb N$. Let $K_\ast:=\dset{v_n}{n \in \mathbb N}\cup \dset{w_n}{n \in \mathbb N} \cup \{u\}$. The set $K_\ast$ is clearly compact. If for any $L\geq 0$ we choose a positive integer $m\geq L$ and the corresponding points $v_m,w_m$, we obtain $\norm{f(v_m)-f(w_m)}>m\norm{v_m-w_m}^{\frac{1}{q}}\geq L\norm{v_m-w_m}^{\frac{1}{q}}$. This implies that $f$ is not H\"older continuous on $K_\ast$, which leads to a contradiction.

We now proceed to the main part of the proof. Let $\varepsilon>0$ and $x \in bv_1(E)$ be given. Set $u_\ast:=\lim_{n \to \infty}x(n)$; note that this limit exists, because $bv_1(E)\subset c(E)$ when the normed space $E$ is complete. Choose $r>0$ and $L_r \geq 0$ such that $\norm{f(w)-f(u)}\leq L_r\norm{w-u}^{\frac{1}{q}}$ for all $w,u \in \ball_E(u_\ast,r)$. Furthermore, select $N \in \mathbb N$ such that $\norm{x(n)-u_\ast}\leq \frac{1}{2}r$ for all $n \geq N$, and
\[
 \sum_{n=N}^\infty \norm{x(n+1)-x(n)} \leq \frac{\varepsilon^q}{2^{2q+3}(L_r+1)^q}.
\]
Finally, let
\[
 0< \delta \leq \min\biggl\{\frac{r}{2}, \frac{\varepsilon^q}{2^{2q+3}(L_r+1)^q}\biggr\}
\]
be such that
\[
 \norm{f(x(n))-f(w)} \leq \frac{\varepsilon}{2^{2+\frac{2}{q}}N^{\frac{1}{q}}}
\]
for all $n=1,\ldots,N$ and $w \in \ball_E(x(n),\delta)$; the existence of such a $\delta$ is guaranteed by the continuity of $f$ on $E$. 

Now, consider any $y \in \ball_{bv_1}(x,\delta)$. Then, we have $\norm{x(n)-y(n)}\leq \delta$ for every $n \in \mathbb N$. Hence, $\norm{y(n)-u_\ast}\leq r$ for $n \geq N$. Also,
\begin{align*}
&\sum_{n=N}^\infty \norm{y(n+1)-y(n)}\\
&\qquad \leq \sum_{n=N}^\infty \norm[\big]{[y(n+1)-x(n+1)]-[y(n)-x(n)]} + \sum_{n=N}^\infty \norm{x(n+1)-x(n)}\\
&\qquad \leq \norm{x-y}_{bv_1} +  \sum_{n=N}^\infty \norm{x(n+1)-x(n)} \leq \frac{\varepsilon^q}{2^{2q+2}(L_r+1)^q}.
\end{align*}
It follows that
\begin{align*}
& \sum_{n=1}^{N-1}\norm[\Big]{\bigl[f(x(n+1))-f(y(n+1))\bigr]-\bigl[f(x(n))-f(y(n))\bigr]}^q \\
&\qquad \leq 2^{q+1} \sum_{n=1}^{N}\norm{f(x(n))-f(y(n))}^q \leq 2^{q+1} N \cdot \frac{\varepsilon^q}{2^{2q+2}N} = \frac{\varepsilon^q}{2^{q+1}}.
\end{align*}
And,
\begin{align*}
& \sum_{n=N}^{\infty}\norm[\Big]{\bigl[f(x(n+1))-f(y(n+1))\bigr]-\bigl[f(x(n))-f(y(n))\bigr]}^q \\
&\qquad \leq 2^q\sum_{n=N}^{\infty}\norm{f(x(n+1))-f(x(n))}^q + 2^q\sum_{n=N}^{\infty}\norm{f(y(n+1))-f(y(n))}^q\\
&\qquad \leq 2^q L_r^q \sum_{n=N}^{\infty} \norm{x(n+1)-x(n)} + 2^q L_r^q \sum_{n=N}^{\infty} \norm{y(n+1)-y(n)}\\
&\qquad \leq 2^{q} L_r^q \cdot  \frac{\varepsilon^q}{2^{2q+3}(L_r+1)^q} + 2^{q} L_r^q \cdot \frac{\varepsilon^q}{2^{2q+2}(L_r+1)^q} = \frac{\varepsilon^q}{2^{q+3}} + \frac{\varepsilon^q}{2^{q+2}} \leq \frac{\varepsilon^q}{2^{q+1}}.
\end{align*}
Combining the above estimates with the inequality $\norm{f(x(1))-f(y(1))} \leq \frac{1}{4}\varepsilon$, yields $\norm{C_f(x)-C_f(y)}_{bv_q} \leq \varepsilon$. Thus, we conclude that the composition operator $C_f$ is continuous at $x$.

\item Let $1<p\leq q<+\infty$. Since this case follows similarly to the previous one, we will omit some minor details. From the characterization of acting conditions for the composition operator $C_f \colon bv_p(E) \to bv_q(E)$, we know that there exist constants $r>0$ and $L\geq 0$ such that $\norm{f(u)-f(w)} \leq L\norm{u-w}^{\frac{p}{q}}$ for all $u,w \in E$ with $\norm{u-w}\leq r$ (see \cite{BK}*{Theorem 23}). Now, let $\varepsilon>0$ and $x \in bv_p(E)$ be fixed. Also, let $N \in \mathbb N$ be chosen so that
\[
 \sum_{n=N}^\infty \norm{x(n+1)-x(n)}^p \leq \min\biggl\{\frac{r^p}{2^p}, \frac{\varepsilon^q}{2^{2q+p+3}(L+1)^q}\biggr\}.
\]
Since the function $f$ is uniformly continuous on $E$ (and hence continuous), there exists a number
\[
 0<\delta \leq \min\biggl\{\frac{r}{2}, \frac{\varepsilon^{\frac{q}{p}}}{2^{\frac{2q+3}{p}+1}(L+1)^\frac{q}{p}} \biggr\}
\]
such that
\[
 \norm{f(x(n))-f(w)}\leq \frac{\varepsilon}{2^{2+\frac{2}{q}}N^\frac{1}{q}}
\]
for all $n=1,\ldots,N$ and $w \in \ball_E(x(n),\delta)$.

Now, take any $y \in \ball_{bv_p}(x,\delta N^{-1})$. Then, $\norm{y(n+1)-y(n)}\leq r$ for $n \geq N$ and
\begin{align*}
& \sum_{n=N}^\infty \norm{y(n+1)-y(n)}^p\\
&\qquad \leq 2^p\sum_{n=N}^\infty \norm[\big]{[y(n+1)-x(n+1)]-[y(n)-x(n)]}^p + 2^p\sum_{n=N}^\infty \norm{x(n+1)-x(n)}^p\\
&\qquad \leq 2^p\delta^p + 2^p \cdot \frac{\varepsilon^q}{2^{2q+p+3}(L+1)^q}= \frac{\varepsilon^q}{2^{2q+2}(L+1)^q}.
\end{align*}
Furthermore, by the definition of the norm in $bv_p(E)$, we have $\norm{x(1)-y(1)}\leq \delta N^{-1}$ and similarly,
\[
 \norm{x(2)-y(2)}\leq \norm[\big]{[x(2)-y(2)]-[x(1)-y(1)]} + \norm{x(1)-y(1)}\leq 2\delta N^{-1}.
\] 
By continuing this process, we obtain $\norm{x(n)-y(n)}\leq \delta$ for $n=1,\ldots,N$. As a result, we have the following estimate
\begin{align*}
& \sum_{n=1}^{N-1}\norm[\Big]{\bigl[f(x(n+1))-f(y(n+1))\bigr]-\bigl[f(x(n))-f(y(n))\bigr]}^q \\
&\qquad \leq 2^{q+1} \sum_{n=1}^{N}\norm{f(x(n))-f(y(n))}^q \leq \frac{\varepsilon^q}{2^{q+1}}.
\end{align*}
And,
\begin{align*}
& \sum_{n=N}^{\infty}\norm[\Big]{\bigl[f(x(n+1))-f(y(n+1))\bigr]-\bigl[f(x(n))-f(y(n))\bigr]}^q \\
&\qquad \leq 2^q\sum_{n=N}^{\infty}\norm{f(x(n+1))-f(x(n))}^q + 2^q\sum_{n=N}^{\infty}\norm{f(y(n+1))-f(y(n))}^q\\
&\qquad \leq 2^q L^q \sum_{n=N}^{\infty} \norm{x(n+1)-x(n)}^p + 2^q L^q \sum_{n=N}^{\infty} \norm{y(n+1)-y(n)}^p \leq \frac{\varepsilon^q}{2^{q+1}}.
\end{align*}
By combining the above estimates with the inequality $\norm{f(x(1))-f(y(1))} \leq \frac{1}{4}\varepsilon$, we conclude that $\norm{C_f(x)-C_f(y)}_{bv_q} \leq \varepsilon$. This completes the proof that the composition operator $C_f$ is continuous at $x$. 
\end{enumerate}
\end{proof}

We end this subsection with two results concerning the pointwise continuity of composition operators acting from $bv_p(E)$ into $c(E)$ and $l^\infty(E)$. What is particularly striking is the stark difference in behaviour between the cases $p=1$ and $p>1$. In the former case, there are non-trivial continuous composition operators, while in the latter no such operators exist.

\begin{theorem}\label{thm:cont_bv_1_c_l}
Let $E$ be a Banach space and let $X \in \{c(E),l^\infty(E)\}$. Moreover, assume that the composition operator $C_f$ maps $bv_1(E)$ into $X$. Then, $C_f$ is continuous if and only if $f$ is.
\end{theorem}

\begin{proof}
The right implication is obvious -- see Proposition~\ref{prop:cont_Cf_implies_f}. Now, assume that $f$ is continuous. Let $\varepsilon>0$ and $x \in bv_1(E)$ be given. Moreover, set $u_\ast:=\lim_{n \to \infty} x(n)$; this limit exists, because $bv_1(E) \subset c(E)$ when $E$ is a Banach space. Since $f$ is continuous at $u_\ast$, there is $\delta>0$ such that $\norm{f(u)-f(u_\ast)}\leq \frac{1}{2}\varepsilon$, whenever $u \in \ball_E(u_\ast,\delta)$. Choose also an index $N \in \mathbb N$ such that $\norm{x(n)-u_\ast}\leq \frac{1}{2}\delta$ for all $n \geq N$. Next, using the continuity of $f$ for the second time, we find $\eta \in (0,\frac{1}{2}\delta]$ such that $\norm{f(x(n))-f(w)}\leq \varepsilon$ for all $n=1,\ldots,N-1$ and $w \in \ball_E(x(n),\eta)$.

Now, consider any $y \in \ball_{bv_1}(x,\eta)$. Then, $\norm{y(n)-u_\ast}\leq \delta$ for $n \geq N$. Hence,
\[
 \norm{f(x(n))-f(y(n))} \leq \norm{f(x(n))-f(u_\ast)} + \norm{f(u_\ast)-f(y(n))} \leq \varepsilon
\]
for all $n\geq N$. For $n=1,\ldots,N-1$ we also have $\norm{f(x(n))-f(y(n))}\leq \varepsilon$ by the choice of $\eta$. Thus, we conclude that $\norm{C_f(x)-C_f(y)}_\infty \leq \varepsilon$, demonstrating that the composition operator $C_f \colon bv_1(E) \to X$ is continuous.
\end{proof}

It is known that the composition operator $C_f$ maps $bv_1(E)$ into $c(E)$, where $E$ is Banach space, if and only if its generator $f \colon E \to E$ is continuous (see~\cite{BK}*{Theorem~29}). Therefore, in view of Theorem~\ref{thm:cont_bv_1_c_l}, we can immediately deduce the following result.

\begin{corollary}
Let $E$ be a Banach space. Ever composition operator $C_f$ that maps $bv_1(E)$ into $c(E)$ is continuous.
\end{corollary}

The situation is quite different for composition operators acting from $bv_p(E)$ when $p>1$.

\begin{theorem}
Let $p \in (1,+\infty)$ and let $X \in \{c(E),l^\infty(E)\}$. Moreover, assume that the composition operator $C_f$ maps $bv_p(E)$ into $X$. Then, $C_f$ is continuous if and only if $f$ is a constant mapping.
\end{theorem}

\begin{proof}
In view of~\cite{BK}*{Theorem~31}, we only need to consider the case $X:=l^\infty(E)$.

Suppose that there exists a point $u \in E$ such that $f(u)\neq f(0)$. Furthermore, fix an $r \in \mathbb N$ such that $\norm{u} \leq 2^{r}$. And for each $n \in \mathbb N$ set
\[
 a_n:=\biggl[ \frac{n+pr}{p-1}\biggr] + 1,
\]
where $[t]$ denotes the smallest integer larger than or equal to $t$. Define the sequence $x$ as
\[
0,\frac{1}{2^{a_1}}u,\frac{2}{2^{a_1}}u,\ldots,\frac{2^{a_1}-1}{2^{a_1}}u,u,\frac{2^{a_1}-1}{2^{a_1}}u,\ldots,\frac{2}{2^{a_1}}u, \frac{1}{2^{a_1}}u, 0, \frac{1}{2^{a_2}}u, \frac{2}{2^{a_2}}u, \ldots, \frac{2^{a_2}-1}{2^{a_2}}u, u,\ldots. 
\]
Note that
\[
 \sum_{n=1}^\infty \norm{x(n+1)-x(n)}^p \leq 2 \sum_{n=1}^\infty 2^{a_n} \cdot 2^{-p a_n}\norm{u}^p \leq 2\sum_{n=1}^\infty 2^{-n} = 2 < +\infty.
\]
This means that $x \in bv_p(E)$. Now, for each $m \in \mathbb N$ set $s_m:=\sum_{i=1}^m 2^{a_i+1}$ and
\[
  x_m:=\bigl(x(1),x(2),\ldots,x(s_m),0,0,0,\ldots \bigr)
\]
Observe that
\[
 \norm{x-x_m}^p_{bv_p} \leq 2\sum_{n=m+1}^\infty 2^{a_n} \cdot 2^{-pa_n}\norm{u}^p \leq 2\sum_{n=m+1}^\infty 2^{-n} = 2^{-m+1}.
\]
Therefore, $x_m \to x$ in $bv_p(E)$. However, $\norm{C_f(x) - C_f(x_m)}_\infty \geq \norm{f(u)-f(0)}>0$ for any $m \in \mathbb N$. Hence, $C_f$ is not continuous at $x$.
\end{proof}

\subsection{Local uniform continuity}
\label{sec:local_unifom_continuity}

As the title suggests, we now turn our attention to characterizing the local uniform continuity of certain composition operators. Unlike the previous subsection, the forthcoming theorems will involve more complex conditions, particularly in the context of composition operator acting between $bv_p$-spaces.

Let us begin, however, with a simple result.

\begin{theorem}
Let $p \in [1,+\infty)$. Moreover, assume that the composition operator $C_f$ maps $c_0(E)$ into $bv_p(E)$. Then, $C_f$ is locally uniformly continuous if and only if $f$ is a constant map.
\end{theorem}

\begin{proof}
Clearly, it suffices to prove only the necessity part. Assume that $C_f$ is locally uniformly continuous. Then, by Proposition~\ref{prop:cont_Cf_implies_f} (or Theorem~\ref{thm:c_0_bvp_continuous}), the generator $f$ is continuous. Moreover, by \cite{BK}*{Theorem~15}, it is also constant on some closed ball $\ball_E(0,\delta)$ of positive radius. Define
\[
 R:=\sup\dset{r>0}{\text{$f|_{\ball_E(0,r)}$ is constant}}.
\]
If $R=+\infty$, there is noting to prove, as $f$ is constant on $E$. So, let us assume that $R<+\infty$, and let $(r_n)_{n \in \mathbb N}$ be a sequence in
\[
 \dset{r>0}{\text{$f|_{\ball_E(0,r)}$ is constant}}
\]
such that $r_n \to R$. For any $w \in \ball_E(0,R)$, define $w_n:=\frac{r_n}{R}w$. Since $w_n \in \ball_E(0,r_n)$ for all $n \in \mathbb N$, it follows that $f(w_n)=f(0)$. As the map $f$ is continuous and $w_n \to w$ as $n \to +\infty$, we conclude that $f(w)=f(0)$. Hence, $f$ is constant on $\ball_E(0,R)$. This, in turn, implies that for each $n \in \mathbb N$ there exists $u_n \in \ball_E(0,R+\frac{1}{n})\setminus \ball_E(0,R)$ such that $f(u_n)\neq f(0)$. Let $v_n:=\frac{R}{\norm{u_n}}u_n$. It is straightforward to check that $\norm{u_n-v_n}\leq \frac{1}{n}$. Next, choose $k_n \in \mathbb N$ such that $k_n^{\frac{1}{p}}\norm{f(u_n)-f(0)}\geq 1$, and define the sequences $x_n:=(u_n,0,u_n,0,\ldots,u_n,0,0,\ldots)$ and $y_n:=(v_n,0,v_n,0,\ldots,v_n,0,0,\ldots)$, where the last non-zero term appears at the $(2k_n-1)$-th position. Observe that $x_n,y_n \in \ball_{c_0}(0,R+1)$ and  $\norm{x_n-y_n}_\infty = \norm{u_n-v_n}\leq \frac{1}{n}$. However, $\norm{C_f(x_n)-C_f(y_n)}_{bv_p} \geq k_n^{\frac{1}{p}}\norm{f(u_n)-f(0)}\geq 1$. This contradicts the uniform continuity of $C_f$ on $\ball_{c_0}(0,R+1)$, completing the proof.
\end{proof}

Next, we aim to characterize locally uniformly continuous composition operators acting between $l^p(E)$ and $bv_q(E)$. 

\begin{theorem}
Let $p,q \in [1,+\infty)$. Moreover, assume that the composition operator $C_f$ maps $l^p(E)$ into $bv_q(E)$. Then, $C_f$ is locally uniformly continuous if and only if for every $\varepsilon>0$ and $r>0$ there exist $\delta>0$ and constants $L_1,L_2,L_3\geq 0$ satisfying $(L_1^q+L_2^q)r^p\leq \varepsilon^q$ such that for all $u,w \in \ball_E(0,r)$ we have
\[
 \norm{f(u)-f(w)}\leq L_1\norm{u}^{\frac{p}{q}}+L_2\norm{w}^{\frac{p}{q}} + L_3\norm{u-w}^{\frac{p}{q}},
\]
whenever $\norm{u-w}\leq \delta$. 
\end{theorem}

\begin{proof}
We begin by proving the local uniform continuity of the composition operator. Fix $\varepsilon>0$ and $r>0$. Choose the numbers $\delta>0$ and $L_1,L_2,L_3 \geq 0$ according to the growth condition of $f$, applied to $r$ and $12^{-(1+\frac{1}{q})}\varepsilon$ in place of $\varepsilon$. Note that, keeping the constant $L_3$ fixed, we may assume that 
\[
 0 < \delta^p < \frac{12^{-(q+1)}\varepsilon^q}{L_3^q + 1}.
\] 
Now, let $x,y \in \ball_{l^p}(0,r)$ satisfy $\norm{x-y}_{l^p}\leq \delta$. Then, for all $n \in \mathbb N$, we have $x(n),y(n) \in \ball_E(0,r)$ and $\norm{x(n)-y(n)}\leq \delta$. Consequently,
\begin{align*}
& \sum_{n=1}^\infty \norm[\Big]{\bigl[f(x(n+1))-f(y(n+1))\bigr]-\bigl[f(x(n))-f(y(n))\bigr]}^q\\
& \qquad \leq  2^{q+1}\sum_{n=1}^\infty \norm[\big]{f(x(n))-f(y(n))}^q\\
& \qquad \leq 6^{q+1}L_1^q\sum_{n=1}^\infty \norm{x(n)}^p +  6^{q+1}L_2^q\sum_{n=1}^\infty \norm{y(n)}^p + 6^{q+1}L_3^q\sum_{n=1}^\infty \norm{x(n)-y(n)}^p\\
& \qquad \leq  6^{q+1}(L_1^q + L_2^q)r^p + 6^{q+1}L_3^q\delta^p \leq 2^{-q}\varepsilon^q. 
\end{align*}
Similarly, we obtain $\norm{f(x(1))-f(y(1))}\leq \frac{1}{2}\varepsilon$. Hence, $\norm{C_f(x)-C_f(y)}_{bv_q}\leq \varepsilon$, which proves that the composition operator $C_f$ is locally uniformly continuous.

Now, we proceed to prove the opposite implication by contradiction. Suppose that the composition operator $C_f$ is locally uniformly continuous, but there exist numbers $\varepsilon>0$, $r>0$ and sequences $(u_n)_{n \in \mathbb N}$ and $(w_n)_{n \in \mathbb N}$ in $\ball_E(0,r)$ such that $\norm{u_n-w_n}\leq 4^{-\frac{1}{p}}n^{-\frac{2q}{p}}$ and
\[
 \norm{f(u_n)-f(w_n)} > 2^{-\frac{1}{q}}\varepsilon r^{-\frac{p}{q}}\norm{u_n}^{\frac{p}{q}} + 2^{-\frac{1}{q}}\varepsilon r^{-\frac{p}{q}}\norm{w_n}^{\frac{p}{q}} + n^2\norm{u_n-w_n}^{\frac{p}{q}}. 
\]
Without loss of generality, we may assume that $\varepsilon \leq \frac{1}{4}$. For each $n \in \mathbb N$ let us set
\[
 m_n:=\bigl[\bigl(\tfrac{1}{2}\varepsilon^q r^{-p}\norm{u_n}^p + \tfrac{1}{2}\varepsilon^q r^{-p}\norm{w_n}^p + n^{2q}\norm{u_n-w_n}^p\bigr)^{-1} \bigr],
\]
where $[t]$ denotes the smallest integer greater than or equal to $t$. Note that 
\[
 0<\tfrac{1}{2}\varepsilon^q r^{-p}\norm{u_n}^p + \tfrac{1}{2}\varepsilon^q r^{-p}\norm{w_n}^p + n^{2q}\norm{u_n-w_n}^p \leq \varepsilon^q + \tfrac{1}{4} \leq \tfrac{1}{2}
\]
for all $n \in \mathbb N$. Hence, the numbers $m_n$ are well-defined and satisfy $m_n \geq 2$. Next, define two sequences 
\begin{align*}
 x_n&:=(0,u_n,0,u_n,0,\ldots,0,u_n,0,0,\ldots),\\
 y_n&:=(0,w_n,0,w_n,0,\ldots,0,w_n,0,0,\ldots),
\end{align*}
where the last element $u_n$ or $w_n$ appears at the $2m_n$-th position. Then, we have $ \norm{x_n}_{l^p}^p = m_n \norm{u_n}^p \leq 2\varepsilon^{-q} r^p$, and similarly, $\norm{y_n}_{l^p}^p \leq 2\varepsilon^{-q}r^p$. Thus, for every $n \in \mathbb N$, both $x_n$ and $y_n$ belong to the closed ball $\ball_{l^p}(0,R)$, where $R:=2^{\frac{1}{p}}\varepsilon^{-\frac{q}{p}}r$. Furthermore, $\norm{x_n -y_n}^p_{l^p} = m_n\norm{u_n-w_n}^p \leq n^{-2q}$. Therefore, $\norm{x_n-y_n}_{l^p} \to 0$ as $n \to +\infty$. However, for each $n \in \mathbb N$ we have
\begin{align*}
   &\sum_{k=1}^\infty \norm[\Big]{ \bigl[f(x_n(k+1))-f(y_n(k+1))\bigr] - \bigl[f(x_n(k))-f(y_n(k))\bigr]}^q\\
	&\qquad \geq m_n \norm[\big]{f(u_n)-f(w_n)}^q\\
	&\qquad \geq \frac{\bigl(2^{-\frac{1}{q}}\varepsilon r^{-\frac{p}{q}}\norm{u_n}^{\frac{p}{q}} + 2^{-\frac{1}{q}}\varepsilon r^{-\frac{p}{q}}\norm{w_n}^{\frac{p}{q}} + n^2\norm{u_n-w_n}^{\frac{p}{q}}\bigr)^q}{\varepsilon^q r^{-p}\norm{u_n}^p + \varepsilon^q r^{-p}\norm{w_n}^p + 2n^{2q}\norm{u_n-w_n}^p}\\
	& \qquad \geq \frac{\tfrac{1}{2}\varepsilon^q r^{-p}\norm{u_n}^p + \tfrac{1}{2}\varepsilon^q r^{-p}\norm{w_n}^{p} + n^{2q}\norm{u_n-w_n}^{p}}{\varepsilon^q r^{-p}\norm{u_n}^p + \varepsilon^q r^{-p}\norm{w_n}^p + 2n^{2q}\norm{u_n-w_n}^p} \geq \frac{1}{2}. 
\end{align*}
This implies that $\norm{C_f(x_n)-C_f(y_n)}_{bv_q}\geq \frac{1}{2}$ for all $n \in \mathbb N$. Therefore, the composition operator $C_f$ is not locally uniformly continuous. This completes the proof.
\end{proof}

We can now move on to discussing the local uniform continuity of composition operators acting between $bv_1(E)$ and $bv_p(E)$. Unfortunately, in this case, we only have sufficient conditions. Furthermore, in the general setting of an arbitrary Banach space $E$, these conditions appear to be far from necessary.

Before we proceed with proving our first result, we need to introduce some notation. If the normed space $E$ is complete, we denote by $\mathcal L(E)$ the Banach space of all linear and continuous operators $T \colon E \to E$, endowed with the operator norm $\norm{T}_{\mathcal L}:=\sup_{\norm{x}\leq 1}\norm{T(x)}$.

\begin{theorem}\label{thm:LUC_bv_1_bv_p}
Let $E$ be a Banach space. If the generator $f \colon E \to E$ is  Fr\'echet differentiable with the derivative $f' \colon E \to \mathcal{L}(E)$ being locally uniformly continuous and locally bounded, then the composition operator $C_f$ maps $bv_1(E)$ into itself and is locally uniformly continuous. 
\end{theorem}

\begin{proof}
The argument will be loosely based on the approach used in the proof of~\cite{DN}*{Theorem~6.63}, which addresses the local uniform continuity of superposition operators in spaces of functions with bounded Wiener variation. A key tool in our reasoning is a Lagrange-type formula for continuously differentiable maps $F \colon B_E(0,\rho) \to E$, where $\rho>0$. It says that for any points $u,w \in B_E(0,\rho)$ we have
\begin{equation}\label{eq:mean_value}
 F(u)-F(w)=\int_0^1 F'(w+t(u-w))(u-w) \textup dt;
\end{equation}
here, the abstract integral is understood in the Riemann sense (see~\cite{DN}*{Lemma~5.40}).

Since the derivative $f'$ is locally bounded, the mean value theorem for vector-valued  functions (see, for example,~\cite{DN}*{Theorem~5.3}) ensures that $f$ satisfies a Lipschitz condition on bounded subsets of $E$. In particular, by~\cite{BK}*{Theorem~20}, the composition operator $C_f$ maps $bv_1(E)$ into itself. 

Now, fix $r>0$ and $\varepsilon>0$. Let $R\geq 1$ be a constant such that $\norm{f'(u)}_{\mathcal L}\leq R$ for all $u \in \ball_E(0,r)$. Choose $0<\delta \leq \frac{1}{4}R^{-1}\varepsilon$ so that for any $u,w \in \ball_E(0,r)$ satisfying $\norm{u-w}\leq \delta$, we have
$\norm{f'(u)-f'(w)}_{\mathcal L}\leq \frac{1}{4}r^{-1}\varepsilon$ and $\norm{f(u)-f(w)}\leq \frac{1}{2}\varepsilon$. Next, consider any two points $x,y \in \ball_{bv_1}(0,r)$ with $\norm{x-y}_{bv_1}\leq \delta$. Then, for every $n \in \mathbb N$ we have $\norm{x(n)-y(n)}\leq \delta$ and $x(n),y(n) \in \ball_E(0,r)$. Furthermore, for each $n \in \mathbb N$, we arrive at the identity
\begin{align*}
 &\bigl[f(x(n+1))-f(y(n+1))\bigr] - \bigl[f(x(n))-f(y(n))\bigr]\\
 &\qquad = \bigl[f(x(n+1))-f(x(n))\bigr] - \bigl[f(y(n+1))-f(y(n))\bigr]\\
 &\qquad = \int_0^1 f'\bigl(x(n)+t(x(n+1)-x(n))\bigr)(x(n+1)-x(n))\textup dt\\
 &\qquad \qquad - \int_0^1 f'\bigl(y(n)+t(y(n+1)-y(n))\bigr)(y(n+1)-y(n))\textup dt\\
 &\qquad = \int_0^1 \bigl[f'\bigl(x(n)+t(x(n+1)-x(n))\bigr)\\
 &\hspace{4cm} - f'\bigl(y(n)+t(y(n+1)-y(n))\bigr)\bigr](x(n+1)-x(n))\textup dt\\
 &\qquad \qquad + \int_0^1 f'\bigl(y(n)+t(y(n+1)-y(n))\bigr)\bigl[(x(n+1)-x(n))-(y(n+1)-y(n))\bigr]\textup dt.
\end{align*}
Since
\begin{align*}
 &\norm[\big]{x(n)+t(x(n+1)-x(n)) - y(n)-t(y(n+1)-y(n))}\\
 &\hspace{3cm}\leq (1-t) \norm[\big]{x(n)-y(n)} + t\norm{x(n+1)-y(n+1)} \leq (1-t)\delta + t\delta =\delta
\end{align*}
for  $t \in [0,1]$, we obtain
\begin{align*}
&\norm[\Big]{\bigl[f(x(n+1))-f(y(n+1))\bigr] - \bigl[f(x(n))-f(y(n))\bigr]}\\
&\qquad \leq \tfrac{1}{4}r^{-1}\varepsilon \norm{x(n+1)-x(n)} + R\norm[\big]{(x(n+1)-x(n))-(y(n+1)-y(n))}.
\end{align*}
Summing over all $n$, we get
\begin{align*}
 &\sum_{n=1}^\infty \norm[\Big]{\bigl[f(x(n+1))-f(y(n+1))\bigr] - \bigl[f(x(n))-f(y(n))\bigr]}\\
 &\qquad \leq \frac{1}{4} r^{-1}\varepsilon \sum_{n=1}^\infty \norm{x(n+1)-x(n)} + R\sum_{n=1}^\infty\norm[\big]{(x(n+1)-y(n+1)) - (x(n)-y(n))}\\
 &\qquad \leq \frac{1}{4} r^{-1}\varepsilon \norm{x}_{bv_1} + R\norm{x-y}_{bv_1} \leq \frac{1}{4}\varepsilon + \frac{1}{4}\varepsilon = \frac{1}{2}\varepsilon.
\end{align*}
Moreover, we also have $\norm{f(x(1))-f(y(1))} \leq \frac{1}{2}\varepsilon$. Consequently, $\norm{C_f(x)-C_f(y)}_{bv_1} \leq \varepsilon$. This shows that the composition operator $C_f$ is uniformly continuous on the ball $\ball_{bv_1}(0,r)$, completing the proof.
\end{proof}

It turns out that, in the special case where $E:=\mathbb R$, the above result can be improved to ensure that the target space is $bv_p$ for any $p \in [1,+\infty)$. Naturally, this improvement comes at a cost: the required estimates will need to be more delicate and will rely on the following technical lemma.

\begin{lemma}\label{lem:technical_estimates}
Let $\alpha \in (0,1]$.
\begin{enumerate}[label=\textup{(\alph*)}]
 \item If $a,b \geq 0$, then $(a+b)^\alpha \leq a^\alpha + b^\alpha$.
 \item If $0\leq b \leq a$, then $a^{\alpha+1} - b^{\alpha+1} \leq (a-b)^{\alpha+1} + (\alpha+1)b^{\alpha}(a-b)$.
 \item If $a,b,c,d \in \mathbb R$ with $a\neq b$ and $c\neq d$, then
\[
 \abs[\bigg]{ \frac{\abs{a-b}^{1+\alpha}}{a-b} - \frac{\abs{c-d}^{1+\alpha}}{c-d}} \leq 8\abs{a-b-c+d}^{\alpha}.
\]
\end{enumerate}
\end{lemma}

\begin{proof}
The proofs of parts~(a) and~(b) are well-known and standard, so we will omit them. We now proceed to the proof of part~(c). Without loss of generality, we may assume that $\abs{a-b}\geq \abs{c-d}>0$. Using part~(b) we have
\begin{align*}
 &  \abs[\bigg]{ \frac{\abs{a-b}^{1+\alpha}}{a-b} - \frac{\abs{c-d}^{1+\alpha}}{c-d}}\\
 &\qquad = \frac{ \abs[\big]{ (a-b)\abs{c-d}^{1+\alpha} - (c-d)\abs{a-b}^{1+\alpha}}}{\abs{a-b}\abs{c-d}}\\
 & \qquad = \frac{ \abs[\big]{ (a-b-c+d)\abs{c-d}^{1+\alpha} + (c-d)\abs{c-d}^{1+\alpha} - (c-d)\abs{a-b}^{1+\alpha}}}{\abs{a-b}\abs{c-d}}\\
 & \qquad \leq \frac{\abs{a-b-c+d}\abs{c-d}^{\alpha}}{\abs{a-b}} + \frac{\abs{a-b}^{1+\alpha} - \abs{c-d}^{1+\alpha}}{\abs{a-b}}\\
 & \qquad \leq \frac{\abs{a-b-c+d}\abs{c-d}^{\alpha}}{\abs{a-b}} + \frac{ \abs{a-b-c+d}^{1+\alpha} + (1+\alpha)\abs{c-d}^{\alpha}\abs{a-b-c+d}}{\abs{a-b}}.
\end{align*}
Now, by the triangle inequality and our assumption, we get $\abs{a-b-c+d} \leq 2\abs{a-b}$. Therefore, we have
\begin{align*}
 \abs{a-b-c+d} &= \abs{a-b-c+d}^{\alpha} \abs{a-b-c+d}^{1-\alpha}\\
 & \leq 2^{1-\alpha}\abs{a-b-c+d}^{\alpha} \abs{a-b}^{1-\alpha}\\
 & \leq 2\abs{a-b-c+d}^{\alpha} \abs{a-b}^{1-\alpha}. 
\end{align*}
Consequently,
\begin{align*}
 &  \abs[\bigg]{ \frac{\abs{a-b}^{1+\alpha}}{a-b} - \frac{\abs{c-d}^{1+\alpha}}{c-d}}\\
& \qquad \leq \frac{2\abs{a-b-c+d}^{\alpha}\abs{a-b}}{\abs{a-b}} + \frac{2\abs{a-b-c+d}^{\alpha}\abs{a-b} + 4\abs{a-b-c+d}^{\alpha}\abs{a-b}}{\abs{a-b}}\\
& \qquad = 8\abs{a-b-c+d}^{\alpha}.
\end{align*}
Thus, we have established the desired inequality, completing the proof.
\end{proof}

\begin{theorem}\label{thm:LUC_bv_1_bv_p_R}
Let $p \in [1,+\infty)$. If $f\colon \mathbb R \to \mathbb R$ is quasi continuously differentiable with exponent $\frac{1}{p}$, then the composition operator $C_f$ maps $bv_1(\mathbb R)$ into $bv_p(\mathbb R)$ and is locally uniformly continuous.    
\end{theorem}

\begin{proof}
By Proposition~\ref{prop:qc1_vs_lip_loc}, the function $f$ is H\"older continuous on bounded subsets of $\mathbb R$ with exponent $\frac{1}{p}$. In particular, in view of~\cite{BK}*{Theorem~20}, the composition operator $C_f$ maps $bv_1(\mathbb R)$ into $bv_p(\mathbb R)$. 

Now, fix $\varepsilon>0$ and $r>0$. Let $L_r \geq 0$ denote the H\"older constant of the function $f$ when restricted to the interval $[-r,r]$. For simplicity, also set
\[
 \gamma:=\frac{\varepsilon}{L_r+2(r+8L_r^p)^{\frac1p}}.
\]
Next, for $u\neq w$ define 
\[
   \Phi(u,w):=\dfrac{(f(u)-f(w))(u-w)}{\abs{u-w}^{1+\frac{1}{p}}}.
\] 
The quasi continuous differentiability of $f$ ensures that $\Phi$ is uniformly continuous on the set $\Delta_r:=([-r,r]\times [-r,r])\setminus \dset{(x,x)}{x \in \mathbb R}$. Therefore, there exists a constant $0<\delta\leq \gamma^p$ such that $\abs{\Phi(a,b)-\Phi(c,d)}\leq \gamma$ for all $a,b,c,d \in [-r,r]$ with $a\neq b$, $c\neq d$, $\abs{a-c}\leq\delta$ and $\abs{b-d}\leq\delta$.

Let $x,y$ be two arbitrary points in $\ball_{bv_1}(0,r)$, and assume that $\norm{x-y}_{bv_1}\leq \delta$. Then, in particular, we have $\abs{x(n)-y(n)}\leq \delta$ and $x(n),y(n)\in [-r,r]$ for every $n \in \mathbb N$. Furthermore, since $f$ is H\"older continuous on the interval $[-r,r]$, we obtain $\abs{f(x(1))-f(y(1))} \leq L_r\abs{x(1)-y(1)}^{\frac{1}{p}}\leq L_r\delta^\frac{1}{p}\leq L_r\gamma$. 

Next, we seek to estimate the quantity
\[
 \abs[\big]{f(x(n+1))-f(y(n+1))-f(x(n))+f(y(n))}^p
\]
for every $n \in \mathbb N$. To proceed, let us fix a positive integer $n$ and consider the following four cases.
\begin{enumerate}[label=\emph{Case \arabic*:}, wide, labelindent=0pt]
 \item assume that $x(n+1)\neq x(n)$ and $y(n+1)\neq y(n)$. Using Lemma~\ref{lem:technical_estimates}, we obtain
\begin{align*}
&\abs[\big]{f(x(n+1))-f(y(n+1))-f(x(n))+f(y(n))}\\[1mm]
&\quad =\abs[\Bigg]{\frac{\bigl(f(x(n+1))-f(x(n))\bigr)\bigl(x(n+1)-x(n)\bigr)}{\abs{x(n+1)-x(n)}^{1+\frac1p}}\cdot\frac{\abs{x(n+1)-x(n)}^{1+\frac1p}}{x(n+1)-x(n)}\\[2mm]
&\hspace{4cm} - \frac{\bigl(f(y)(n+1)-f(y)(n)\bigr)\bigl(y(n+1)-y(n)\bigr)}{\abs{y(n+1)-y(n)}^{1+\frac1p}}\cdot\frac{\abs{y(n+1)-y(n)}^{1+\frac1p}}{y(n+1)-y(n)}}\\[2mm]
&\quad = \abs[\Bigg]{\Phi(x(n+1),x(n)) \cdot \frac{\abs{x(n+1)-x(n)}^{1+\frac1p}}{x(n+1)-x(n)}-\Phi(y(n+1),y(n)) \cdot \frac{\abs{y(n+1)-y(n)}^{1+\frac1p}}{y(n+1)-y(n)}}\\[2mm]
&\quad \leq \abs[\big]{\Phi(x(n+1),x(n))-\Phi(y(n+1),y(n))} \cdot \abs[\big]{x(n+1)-x(n)}^{\frac1p}\\[2mm]
&\hspace{4cm} +\abs[\big]{\Phi(y(n+1),y(n))}\cdot \abs[\Bigg]{\frac{\abs{x(n+1)-x(n)}^{1+\frac1p}}{x(n+1)-x(n)}-\frac{\abs{y(n+1)-y(n)}^{1+\frac1p}}{y(n+1)-y(n)}}\\
&\quad \leq \abs[\big]{\Phi(x(n+1),x(n))-\Phi(y(n+1),y(n))} \cdot \abs[\big]{x(n+1)-x(n)}^{\frac1p}\\
&\hspace{4cm} + 8L_r\abs[\big]{x(n+1)-x(n)-y(n+1)+y(n)}^{\frac1p}.
\end{align*}
Thus, we have
\begin{align*}
& \abs[\big]{f(x(n+1))-f(y(n+1))-f(x(n))+f(y(n))}^p\\ 
&\quad \leq 2^p \abs[\big]{\Phi(x(n+1),x(n))-\Phi(y(n+1),y(n))}^p \cdot \abs[\big]{x(n+1)-x(n)}\\
&\qquad +2^{p+3}L_r^p\abs[\big]{x(n+1)-x(n)-y(n+1)+y(n)}\\
&\quad \leq 2^p \gamma^p\abs[\big]{x(n+1)-x(n)} +2^{p+3}L_r^p\abs[\big]{x(n+1)-x(n)-y(n+1)+y(n)}.
\end{align*}

\item assume that $x(n+1)= x(n)$ and $y(n+1)\neq y(n)$. In this case, we have
\begin{align*}
& \abs[\big]{f(x(n+1))-f(y(n+1))-f(x(n))+f(y(n))}^p\\ 
&\quad \leq L_r^p\abs{y(n+1)-y(n)}\\
&\quad \leq 2^p \gamma^p\abs{x(n+1)-x(n)} +2^{p+3}L_r^p\abs[\big]{x(n+1)-x(n)-y(n+1)+y(n)}.
\end{align*}

\item assume that $x(n+1)\neq x(n)$ and $y(n+1)= y(n)$. Similarly to the previous case, we obtain
\begin{align*}
& \abs[\big]{f(x(n+1))-f(y(n+1))-f(x(n))+f(y(n))}^p\\ 
&\quad \leq 2^p \gamma^p\abs{x(n+1)-x(n)} +2^{p+3}L_r^p\abs[\big]{x(n+1)-x(n)-y(n+1)+y(n)}.
\end{align*}

\item assume that $x(n+1)=x(n)$ and $y(n+1)=y(n)$. In this case, trivially, 
\[
 \abs[\big]{f(x(n+1))-f(y(n+1))-f(x(n))+f(y(n))}^p=0.
\]
Thus, we have
\begin{align*}
& \abs[\big]{f(x(n+1))-f(y(n+1))-f(x(n))+f(y(n))}^p\\ 
&\quad \leq 2^p \gamma^p\abs{x(n+1)-x(n)} +2^{p+3}L_r^p\abs[\big]{x(n+1)-x(n)-y(n+1)+y(n)}.
\end{align*}
\end{enumerate}
By combining all the cases above and summing over all $n \in \mathbb N$, we obtain the following estimate
\begin{align*}
& \sum_{n=1}^{\infty}\abs[\big]{f(x(n+1))-f(y(n+1))-f(x(n))+f(y(n))}^p\\ 
&\quad \leq 2^p \gamma^p \sum_{n=1}^{\infty}\abs{x(n+1)-x(n)}+2^{p+3}L_r^p\sum_{n=1}^{\infty}\abs[\big]{x(n+1)-x(n)-y(n+1)+y(n)}\\
&\quad \leq 2^p \gamma^p\norm{x}_{bv_1}+2^{p+3}L_r^p\norm{x-y}_{bv_1}.
\end{align*}
Consequently,
\begin{align*}
& \norm{C_f(x)-C_f(y)}_{bv_p}\leq L_r\gamma+(2^p\gamma^pr+2^{p+3}L_r^p\gamma^p)^{\frac1p} = \gamma \bigl[L_r+2(r+8L_r^p)^{\frac1p}\bigr]=\varepsilon.
\end{align*}
Thus, the composition operator $C_f$ is locally uniformly continuous. 
\end{proof}

\begin{remark}\label{rem:LUC_approach}
The technique developed in the proof of Theorem~\ref{thm:LUC_bv_1_bv_p_R} appears to be also applicable in providing sufficient conditions for the local uniform continuity of composition operators between $BV_p$-spaces of functions of bounded Wiener variation.  
\end{remark}

Combining Theorem~\ref{thm:LUC_bv_1_bv_p_R} with the results from Section~2.2, we immediately obtain the following corollary.

\begin{corollary}\label{cor:LUC_bv_1_bv_p_R}
Let $p \in [1,+\infty)$. The composition operator $C_f$ maps $bv_1(\mathbb R)$ into $bv_p(\mathbb R)$ and is locally uniformly continuous, if the function $f \colon \mathbb R \to \mathbb R$ satisfies one of the following conditions\textup:
\begin{enumerate}[label=\textup{(\alph*)}]
  \item $f$ is continuously differentiable,
	\item $f$ is locally little H\"older continuous with exponent $\alpha$ such that $\frac{1}{p}\leq\alpha\leq 1$,
	\item $f$ is H\"older continuous on bounded subsets of $\mathbb R$ with exponent $\alpha$ such that $\frac{1}{p}<\alpha\leq 1$.
\end{enumerate}
\end{corollary}

As shown in Section~2.2, for the exponent $\alpha=1$, the classes $C^1(\mathbb R)$ and $\qC^{1,\alpha}(\mathbb R)$ coincide. This naturally leads to the question: in this case, does the converse of Theorem~\ref{thm:LUC_bv_1_bv_p_R} hold, providing a complete characterization of the local uniform continuity of composite operators in the space $bv_1(\mathbb R)$? The answer is affirmative, and the corresponding result is presented below.

\begin{theorem}\label{thm:LUC_C1}
A composition operator $C_f$ mapping $bv_1(\mathbb R)$ into itself is locally uniformly continuous if and only if its generator $f$ is continuously differentiable.  
\end{theorem}

\begin{proof}
Let $C_f$ be a composition operator mapping $bv_1(\mathbb R)$ into itself. If the function $f$ is continuously differentiable, then by Corollary~\ref{cor:LUC_bv_1_bv_p_R}, it follows that $C_f$ is locally uniformly continuous. 

Now, we prove the converse implication. Our approach is inspired by~\cite{reinwand}*{Lemma~5.1.18 and Theorem~5.1.23}. We also rely on the result concerning absolutely continuous functions mentioned in the proof of Proposition~\ref{prop:qc1_vs_c1}. Assume that $C_f$ is locally uniformly continuous. From the characterization of acting conditions for $C_f \colon bv_1(\mathbb R) \to bv_1(\mathbb R)$ it follows that $f$ is Lipschitz continuous on every compact subset of $\mathbb R$, which in this case is equivalent to being Lipschitz on all bounded subsets (see~\cite{BK}*{Theorem~20}). Consequently, $f$ is absolutely continuous on every non-empty and bounded interval. Let $D_f$ denote the set of points in $\mathbb R$ where $f$ is differentiable.  

Fix $r \geq \frac{1}{16}$. We will show that the restriction of $f'$ to the set $[-r,r] \cap D_f$ is uniformly continuous. This will imply that $f$ is continuously differentiable on the entire interval $[-r,r]$. Since $r$ was chosen arbitrarily, the theorem would then follow (cf. the proof of Proposition~\ref{prop:qc1_vs_c1} or \cite{reinwand}*{Lemma~1.1.27~(a)}).

Given $\varepsilon>0$, let $\delta>0$ be such that $\norm{C_f(x)-C_f(y)}_{bv_1}\leq \varepsilon$ for any points $x,y$ in $\ball_{bv_1}(0,33r)$ satisfying $\norm{x-y}_{bv_1}\leq \delta$. Further, fix $u,w \in [-r,r]\cap D_f$ with $\abs{u-w}\leq \delta$, and let $\rho \in (0,r]$. Choose a positive integer $m$ such that $m \leq \frac{16r}{\rho} \leq 2m$. Define $x_\rho, y_\rho \colon \mathbb N \to \mathbb R$ by the formulas
\begin{align*}
 x_\rho:&=(u,u+\rho,u,u+\rho,\ldots,u,u+\rho,u,u,u,\ldots),\\
 y_\rho:&=(w,w+\rho,w,w+\rho,\ldots,w,w+\rho,w,w,w,\ldots),
\end{align*}
where the last term of the form $u+\rho$ or $w+\rho$ appears at the $2m$-th position. Then, we have 
\begin{align*}
 \sum_{n=1}^\infty \abs{x_\rho(n+1)-x_\rho(n)}=2m\rho \leq 32r.
\end{align*}
Thus, $\norm{x_\rho}_{bv_1} \leq 33r$, and similarly, $\norm{y_\rho}_{bv_1}\leq 33r$. Moreover, it is straightforward to show that $\norm{x_\rho-y_\rho}_{bv_1}=\abs{u-w}\leq \delta$. Therefore, by the local uniform continuity of the composition operator $C_f$, we obtain
\begin{align*}
& \varepsilon \geq \norm{C_f(x_\rho)-C_f(y_\rho)}_{bv_1} \geq 2m\abs[\Big]{[f(u+\rho)-f(u)] - [f(w+\rho)-f(w)]}.
\end{align*}
Consequently,
\[
 \abs[\Bigg]{\frac{f(u+\rho)-f(u)}{\rho} - \frac{f(w+\rho)-f(w)}{\rho}} \leq \varepsilon 
\]
for every $\rho \in (0,r]$. Taking the limit as $\rho \to 0^+$ in the above inequality, we finally obtain $\abs{f'(u)-f'(w)}\leq \varepsilon$. This shows that the derivative $f'$, when restricted to the set $[-r,r]\cap D_f$, is uniformly continuous. Hence, $f$ is continuously differentiable on $\mathbb R$ (cf. the proof of Proposition~\ref{prop:qc1_vs_c1} or~\cite{reinwand}*{Lemma~1.1.27~(a)}). 
\end{proof} 

At this point, we should proceed to discuss the local uniform continuity of composition operators $C_f \colon bv_p(E) \to bv_q(E)$ for $1<p\leq q<+\infty$. However, we are currently unaware of any natural classes of maps, aside from affine-like functions, that would generate such operators (cf. Theorem~\ref{thm:UC_bvp_bvq}). Consequently, we skip this discussion and instead proceed directly to the final result of this section.

\begin{theorem}\label{thm:LUC_bv_1_c}
Let $E$ be a Banach space and let $X \in \{c(E), l^\infty(E)\}$. Moreover, assume that the composition operator $C_f$ maps $bv_1(E)$ into $X$. Then, $C_f$ is locally uniformly continuous if and only if $f$ is.
\end{theorem}

\begin{proof}
In view of Proposition~\ref{prop:cont_Cf_implies_f} it is sufficient to prove only the left implication. Therefore, assume that $f$ is locally uniformly continuous on $E$. Fix $r>0$ and $\varepsilon>0$. By the assumption, there exists $\delta>0$ such that $\norm{f(u)-f(w)}\leq \varepsilon$ for all $u,w \in \ball_E(0,r)$ with $\norm{u-w}\leq \delta$. Now, take any $x,y \in \ball_{bv_1}(0,r)$ such that $\norm{x-y}_{bv_1}\leq \delta$. Since $\norm{x}_{\infty}\leq \norm{x}_{bv_1}$, it follows that $x(n),y(n) \in \ball_E(0,r)$ and $\norm{x(n)-y(n)}\leq \delta$ for every $n \in \mathbb N$. Consequently, $\norm{f(x(n))-f(y(n))}\leq \varepsilon$ for $n \in \mathbb N$. This implies that $\norm{C_f(x)-C_f(y)}_\infty \leq \varepsilon$ and shows that $C_f$ is uniformly continuous on $\ball_{bv_1}(0,r)$. 
\end{proof}

\subsection{Uniform continuity}

Finally, in this brief section, we discus uniformly continuous composition operators. The first result is quite straightforward. Since its proof closely mirrors that of Theorem~\ref{thm:LUC_bv_1_c}, we have decided to omit it. The only necessary modification is to replace any ball in $E$ or $bv_1(E)$ with the entire space.

\begin{theorem}\label{thm:UC_bv_1_c}
Let $E$ be a Banach space and let $X \in \{c(E), l^\infty(E)\}$. Moreover, assume that the composition operator $C_f$ maps $bv_1(E)$ into $X$. Then, $C_f$ is uniformly continuous if and only if $f$ is.
\end{theorem}

In contrast to the previous theorem, characterizing the uniform continuity of composition operators between $l^p(E)$ and $bv_q(E)$  requires more careful analysis. 

\begin{theorem}\label{thm:UC_lp_bvq}
Let $p,q \in [1,+\infty)$. Moreover, assume that the composition operator maps $l^p(E)$ into $bv_q(E)$. Then, $C_f$ is uniformly continuous if and only if $f$ is locally H\"older continuous in the stronger sense with exponent $\frac{p}{q}$.
\end{theorem}

\begin{proof}
First, we demonstrate that the composition operator $C_f$ is uniformly continuous. Fix $\varepsilon>0$, and let $L\geq 0$ and $\delta>0$ be the constants appearing in the definition of the local H\"older condition in the stronger sense with exponent $\frac{p}{q}$. Set
\[
 \eta:=\min\Bigl\{\delta, 2^{-\frac{2q+1}{p}}\varepsilon^{\frac{q}{p}} (L+1)^{-\frac{q}{p}}\Bigr\}.
\]
Now, consider any two points $x,y \in l^p(E)$ such that $\norm{x-y}_{l^p}\leq \eta$. Then, $\norm{x(n)-y(n)}\leq \eta \leq \delta$ for all $n \in \mathbb N$. Consequently, we obtain $\norm{f(x(1))-f(y(1))} \leq L\norm{x(1)-y(1)}^{\frac{p}{q}} \leq L\eta^{\frac{p}{q}}\leq \frac{1}{2}\varepsilon$. Moreover,
\begin{align*}
  & \sum_{n=1}^\infty \norm[\Big]{\bigl[f(x(n+1))-f(y(n+1))\bigr] - \bigl[f(x(n))-f(y(n))\bigr]}^q\\
	&\qquad \leq 2^{q+1}\sum_{n=1}^\infty \norm[\big]{f(x(n))-f(y(n))}^q \leq 2^{q+1}L^q\sum_{n=1}^\infty \norm{x(n)-y(n)}^p \leq 2^{q+1}L^q\eta^p \leq \frac{\varepsilon^q}{2^q}.
\end{align*}
Thus, we conclude that $\norm{C_f(x)-C_f(y)}_{bv_q}\leq \varepsilon$, proving that the composition operator $C_f$ is uniformly continuous.

The proof of the second part follows a similar approach to~\cite{BK}*{Theorem 23},  with some minor yet significant modifications. Assume, for the sake of contradiction, that the generator $f$ is not locally H\"older continuous in the stronger sense with exponent $\frac{p}{q}$. Then, there exist two sequences $(u_n)_{n \in \mathbb N}$ and $(w_n)_{n \in \mathbb N}$ in $E$ such that
\[
 \norm{u_n-w_n}\leq 2^{-\frac{1}{p}}n^{-\frac{2q}{p}} \quad \text{and} \quad \norm{f(u_n)-f(w_n)} > n^{2}\norm{u_n-w_{n}}^{\frac{p}{q}}
\]
for all $n \in \mathbb N$. For each $n \in \mathbb N$ set $m_n:=[n^{-2q}\norm{u_n-w_n}^{-p}]$, where $[t]$ denotes the greatest integer less than or equal to $t$. Next, define the sequences $x_n,y_n \colon \mathbb N \to E$ by the formulas
\begin{align*}
 x_n&:=(0,u_n,0,u_n,\ldots,0, u_n,0,0,0,\ldots),\\
 y_n&:=(0,w_n,0,w_n,\ldots,0, w_n,0,0,0,\ldots),
\end{align*}
where $u_n$ and $w_n$ each appear exactly $m_n$ times. Clearly, $x_n, y_n  \in l^p(E)$. Moreover, 
\begin{align*}
 \norm{x_n-y_n}_{l^p}^p = \sum_{k=1}^\infty \norm{x_n(k)-y_n(k)}^p = m_n \norm{u_n-w_n}^p \leq n^{-2q}.
\end{align*}
Thus, $\norm{x_n-y_n}_{l^p} \to 0$ as $n \to +\infty$. However, for every $k \in \{1,\ldots, 2m_n\}$ we have
\begin{align*}
 &\norm[\Big]{\bigl[f(x_n(k+1))-f(y_n(k+1))\bigr] - \bigl[f(x_n(k))-f(y_n(k))\bigr]}^q=\norm[\big]{f(u_n)-f(w_n)}^q.
\end{align*}
Summing over $k$, we obtain
\begin{align*}
 & \sum_{k=1}^\infty \norm[\Big]{\bigl[f(x_n(k+1))-f(y_n(k+1))\bigr] - \bigl[f(x_n(k))-f(y_n(k))\bigr]}^q\\
 &\qquad \geq \sum_{k=1}^{2m_n}\norm[\Big]{\bigl[f(x_n(k+1))-f(y_n(k+1))\bigr] - \bigl[f(x_n(k))-f(y_n(k))\bigr]}^q\\
 & \qquad = \sum_{k=1}^{2m_n} \norm[\big]{f(u_n)-f(w_n)}^q \geq 2m_n n^{2q}\norm{u_n-w_n}^p \geq 1.
\end{align*}
Therefore, $\norm{C_f(x_n)-C_f(y_n)}_{bv_q}\geq 1$ for all $n \in \mathbb N$. This implies that $C_f$ is not uniformly continuous. The proof is complete. 
\end{proof}

Since for the exponent $\alpha>1$ all the H\"older classes discussed in this paper consist solely of constant maps, from Theorem~\ref{thm:UC_lp_bvq} we can deduce the following corollary.

\begin{corollary}
Let $1\leq q < p < +\infty$. Moreover, assume that the composition operator maps $l^p(E)$ into $bv_q(E)$. Then, $C_f$ is uniformly continuous if and only if $f$ is a constant mapping.
\end{corollary}

In the previous section, we observed that providing a complete description of the local uniform continuity of composition operators in $bv_p$-spaces was challenging, and in some cases, even impossible. However, in the current context, this task will be comparatively easier and, in some sense, surprising. The phrase ``in some sense'' is used because those who have worked with superposition operators in spaces of Lipschitz and absolutely continuous functions, or functions of bounded variation, have already encountered similar degeneracy results (see, for example,~\cite{ABM}*{Sections~5.4, 6.2 and~6.3},~\cite{ABKR}*{Theorem~3.3.6},~\cite{BBK}*{Theorem~6.4.16}, or~\cite{reinwand}*{Theorems~5.1.22, 5.2.28 and~5.2.29}). To the best of our knowledge, all of these results were established for $E:=\mathbb R$. For completeness, we note that the first result of this type for a non-autonomous superposition operator acting in the space of functions of bounded Jordan variation was established by Matkowski and Mi\'s in\cite{MatMis}.

Before we state and prove our result, let us recall one last definition. A map $T \colon E \to E$ acting on a normed space $E$ is called $\mathbb R$-linear, if it is additive and satisfies $T(\lambda u)=\lambda T(u)$ for every $u \in E$ and $\lambda \in \mathbb R$. Since in the paper we consider normed spaces over  either the field of real or complex numbers, the notions of an $\mathbb R$-linear and a linear operator may  not always coincide (the former is clearly the broader concept). One of the simplest examples of an $\mathbb R$-linear operator that is not linear is given by the map $T(z)=\operatorname{re}z$ for $z\in \mathbb C$. (Note that $iT(i)=0$, whereas $T(i^2)=-1$.)  

\begin{theorem}\label{thm:UC_bvp_bvq}
Let $1\leq p \leq q < +\infty$. Moreover, assume that $E$ is a Banach space when $p=1$, and a normed space when $p>1$. Assume also that the composition operator $C_f$ maps $bv_p(E)$ into $bv_q(E)$. Then, the following conditions are equivalent:
\begin{enumerate}[label=\textup{(\roman*)}]
 \item\label{it:a} $C_f$ is Lipschitz continuous,
 \item\label{it:b} $C_f$ is uniformly continuous,
 \item\label{it:c} the map $\varphi(u):=f(u)-f(0)$ is a continuous $\mathbb R$-linear endomorphism of $E$.
\end{enumerate}
\end{theorem}

\begin{proof}
The implication $\ref{it:a}~\Rightarrow~\ref{it:b}$ is immediate. Thus, we proceed directly to the proof of the implication $\ref{it:b}~\Rightarrow~\ref{it:c}$. Assume that $C_f$ is uniformly continuous, and let $\omega_{C_f}$ denote its minimal modulus of continuity. That is, for every $\delta>0$, we set
\[
 \omega_{C_f}(\delta):=\sup\dset[\big]{\norm{C_f(x)-C_f(y)}_{bv_q}}{\text{$x,y \in bv_p(E)$ and $\norm{x-y}_{bv_p}\leq \delta$}}.
\] 
It is easy to verify that for every $\delta>0$ the value $\omega_{C_f}(\delta)$ is finite (cf.~\cite{DL}*{Section~2.6}). Now, fix $u,w \in E$, and for each $n \in \mathbb N$ consider two sequences
\begin{align*}
 x_n&:=(\tfrac{1}{2}(u+w), w, \tfrac{1}{2}(u+w), w,\ldots,\tfrac{1}{2}(u+w), w,w,w,\ldots),\\
 y_n&:=(u,\tfrac{1}{2}(u+w), u, \tfrac{1}{2}(u+w),\ldots,u, \tfrac{1}{2}(u+w),\tfrac{1}{2}(u+w),\tfrac{1}{2}(u+w),\tfrac{1}{2}(u+w),\ldots),
\end{align*}
where the ``constant'' part begins with the $2n$-th term. Clearly, $x_n,y_n \in bv_p(E)$, and we have $\norm{x_n-y_n}_{bv_p}=\frac{1}{2}\norm{u-w}$. Furthermore, it follows that
\[
 \omega_{C_f}(\tfrac{1}{2}\norm{u-w}) \geq \norm{C_f(x_n)-C_f(y_n)}_{bv_q} \geq n^{\frac{1}{q}}\norm{2f(\tfrac{1}{2}(u+w))- f(u)-f(w)}.
\]
Dividing the above inequality by $n^{\frac{1}{q}}$ and taking the limit as $n \to +\infty$, we deduce that $f(\tfrac{1}{2}(u+w))=\frac{1}{2}f(u)+\frac{1}{2}f(w)$ for any $u,w \in E$. This implies that $\varphi(u+w)=\varphi(u)+\varphi(w)$ for $u,w \in E$, meaning that the map $\varphi$ is additive. Since $\varphi$ is also continuous by Proposition~\ref{prop:cont_Cf_implies_f}, we conclude that $\varphi(\lambda u)=\lambda \varphi(u)$ for $u \in E$ and $\lambda \in \mathbb R$. (To prove this, it suffices to apply the same reasoning used to show the homogeneity of the solution of the Cauchy functional equation -- cf.~\cite{small}*{Section~2.1}.)

Finally, we prove that the implication $\ref{it:c}~\Rightarrow~\ref{it:a}$ holds. Assume that $\varphi$ is continuous and $\mathbb R$-linear. It is straightforward to show, as in the case of linear mappings, that the continuity of $\varphi$ implies the existence of a constant $L\geq 0$ such that $\norm{\varphi(u)}\leq L\norm{u}$ for all $u \in E$. Let $x,y \in bv_p(E)$. First, we have
\[
 \norm{f(x(1))-f(y(1))}=\norm{\varphi(x(1))-\varphi(y(1))}=\norm{\varphi(x(1)-y(1))} \leq L\norm{x(1)-y(1)}.
\]  
For each $n \in \mathbb N$ we also get
\begin{align*}
& \norm[\big]{f(x(n+1))-f(y(n+1))-f(x(n))+f(y(n))}^q\\
&\qquad = \norm[\big]{\varphi(x(n+1))-\varphi(y(n+1))-\varphi(x(n))+\varphi(y(n))}^q\\
&\qquad \leq L^q\norm[\big]{x(n+1)- y(n+1)-x(n)+y(n)}^q.
\end{align*}
Thus, combining these estimates, we obtain $\norm{C_f(x)-C_f(y)}_{bv_q} \leq L\norm{x-y}_{bv_q} \leq L\norm{x-y}_{bv_p}$, where the last inequality follows from the classical inequality between sums (cf.~\cite{HLP}*{Theorem~19, p.~28}). This shows that the composition operator $C_f \colon bv_p(E) \to bv_q(E)$ is Lipschitz continuous.
\end{proof}

By applying Theorem~\ref{thm:UC_bvp_bvq} to a one-dimensional normed space, we obtain the following result. 

\begin{corollary}\label{cor:UC_p_q_R}
Let $1\leq p \leq q < +\infty$. The composition operator $C_f$ maps $bv_p(\mathbb R)$ into $bv_q(\mathbb R)$ and is uniformly continuous (or, equivalently, Lipschitz continuous) if and only if its generator $f\colon \mathbb R \to \mathbb R$ is of the form $f(u)=au+b$ for some $a,b \in \mathbb R$.
\end{corollary}

\subsection{Continuity: summary}
Let $E$ be a normed space over the real or complex field. If $E$ is a Banach space, we denote it by $\hat E$ to emphasise its completeness. For two sequence spaces $X$ and $Y$, we define the following classes of maps $f \colon E \to E$ that generate composition operators $C_f \colon X \to Y$ with certain continuity properties: 
\begin{align*}
 C(X,Y)&:=\dset{f\colon E \to E}{\text{$C_f$ maps $X$ into $Y$ and is pointwise continuous}},\\
 UC_{\text{loc}}(X,Y)&:=\dset{f\colon E \to E}{\text{$C_f$ maps $X$ into $Y$ and is locally uniformly continuous}},\\
UC(X,Y)&:=\dset{f\colon E \to E}{\text{$C_f$ maps $X$ into $Y$ and is uniformly continuous}}.
\end{align*} 
To fully characterize these classes, we will also incorporate results on acting conditions from~\cite{BK}*{Section~3}.

\clearpage

\begin{small}
\begin{table}[hbt!]
\begin{center}
\begin{tabular}{@{}ccc@{}}
\toprule
\parbox[c]{4.2cm}{\center $f \in C(c_0(E),bv_p(E))$\\ for $1 \leq p < +\infty$} & $\Leftrightarrow$ & \parbox[c]{10.5cm}{\center $f$ is continuous and there exists $\delta>0$\\ such that the restriction $f|_{\ball_(0,\delta)}$ is constant}\\[3mm] \midrule 
\parbox[c]{4.2cm}{\center $f \in C(c(E),bv_p(E))$\\ for $1\leq p < +\infty$} & $\Leftrightarrow$ & $f$ is constant\\[3mm] \midrule 
\parbox[c]{4.2cm}{\center $f \in C(l^\infty(E),bv_p(E))$\\ for $1\leq p < +\infty$} & $\Leftrightarrow$ & $f$ is constant\\[3mm] \midrule 
\parbox[c]{4.2cm}{\center $f \in C(l^p(E),bv_q(E))$\\ for $1\leq p,q < +\infty$} & $\Leftrightarrow$ & \parbox[c]{10.5cm}{\center $f$ is continuous and there exist $\delta>0$ and $L\geq 0$ such that $\norm{f(u)-f(w)}\leq L\norm{u}^{\frac{p}{q}}+L\norm{w}^{\frac{p}{q}}$ for all $u,w \in \ball_E(0,\delta)$} \\[3mm] \midrule 
\parbox[c]{4.2cm}{\center $f \in C(bv_p(E),bv_q(E))$\\ for $1\leq q < p <+\infty$} & $\Leftrightarrow$ & $f$ is constant \\[3mm] \midrule 
\parbox[c]{4.2cm}{\center $f \in C(bv_1(\hat E),bv_p(\hat E))$\\ for $1\leq p <+\infty$} & $\Leftrightarrow$ & \parbox[c]{10.5cm}{\center for every compact subset $K$ of $\hat E$ there exists $L_K\geq 0$ such that $\norm{f(u)-f(w)}\leq L_K\norm{u-w}^{\frac{1}{p}}$ for all $u,w \in K$} \\[4mm] \midrule 
\parbox[c]{4.2cm}{\center $f \in C(bv_p(E),bv_q(E))$\\ for $1< p \leq q<+\infty$} & $\Leftrightarrow$ & \parbox[c]{10.5cm}{\center there exist $\delta>0$ and $L\geq 0$ such that $\norm{f(u)-f(w)}\leq L\norm{u-w}^{\frac{p}{q}}$ for all $u,w \in E$ with $\norm{u-w}\leq \delta$} \\[3.5mm] \midrule 
\parbox[c]{4.2cm}{\center $f \in C(bv_p(E),c_0(E))$\\ for $1 \leq p <+\infty$} & $\Leftrightarrow$ & $f$ is the zero map \\[3mm] \midrule 
\parbox[c]{4.2cm}{\center $f \in C(bv_p(E),l^q(E))$\\ for $1 \leq p,q <+\infty$} & $\Leftrightarrow$ & $f$ is the zero map \\[3mm] \midrule 
$f \in C(bv_1(\hat E),l^\infty(\hat E))$ & $\Leftrightarrow$ & $f$ is continuous on $\hat E$\\ \midrule
\parbox[c]{4.2cm}{\center $f \in C(bv_p(E),l^\infty(E))$\\ for $1<p<+\infty$} & $\Leftrightarrow$ & $f$ is constant\\[3mm] \midrule
$f \in C(bv_1(\hat E),c(\hat E))$ & $\Leftrightarrow$ & $f$ is continuous on $\hat E$ \\ \midrule
\parbox[c]{4.2cm}{\center $f \in C(bv_p(E),c(E))$\\ for $1<p<+\infty$} & $\Leftrightarrow$ & $f$ is constant\\[2mm]
\bottomrule
\end{tabular}
\vspace{1.5mm}
\caption{Pointwise continuity for the composition operator $C_f$.}
\end{center}
\end{table}
\end{small}

\clearpage

\begin{small}
\begin{table}[hbt!]
\begin{center}
\begin{tabular}{@{}ccc@{}}
\toprule
\parbox[c]{4.4cm}{\center $f \in UC_{\text{loc}}(c_0(E),bv_p(E))$\\ for $1 \leq p < +\infty$} & $\Leftrightarrow$ & $f$ is constant\\[3mm] \midrule 
\parbox[c]{4.4cm}{\center $f \in UC_{\text{loc}}(c(E),bv_p(E))$\\ for $1\leq p < +\infty$} & $\Leftrightarrow$ & $f$ is constant\\[3mm] \midrule 
\parbox[c]{4.4cm}{\center $f \in UC_{\text{loc}}(l^\infty(E),bv_p(E))$\\ for $1\leq p < +\infty$} & $\Leftrightarrow$ & $f$ is constant\\[3mm] \midrule 
\parbox[c]{4.4cm}{\center $f \in UC_{\text{loc}}(l^p(E),bv_q(E))$\\ for $1\leq p,q < +\infty$} & $\Leftrightarrow$ & \parbox[c]{12cm}{\center for every $\varepsilon>0$ and $r>0$ there exist $\delta>0$ and constants $L_1,L_2,L_3\geq 0$ satisfying $(L_1^q+L_2^q)r^p\leq \varepsilon^q$ such that for all $u,w \in \ball_E(0,r)$ we have $\norm{f(u)-f(w)}\leq L_1\norm{u}^{\frac{p}{q}}+L_2\norm{w}^{\frac{p}{q}} + L_3\norm{u-w}^{\frac{p}{q}}$,\\ whenever $\norm{u-w}\leq \delta$ 
} \\[4mm] \midrule 
\parbox[c]{4.4cm}{\center $f \in UC_{\text{loc}}(bv_p(E),bv_q(E))$\\ for $1\leq q < p <+\infty$} & $\Leftrightarrow$ & $f$ is constant \\[3mm] \midrule 
\parbox[c]{4.4cm}{\center $f \in UC_{\text{loc}}(bv_1(\hat E),bv_1(\hat E))$} & $\Leftarrow$ & \parbox[c]{10.5cm}{\center $f$ is Fr\'echet differentiable with the derivative $f'$ being locally uniformly continuous and locally bounded} \\[4mm] \midrule 
\parbox[c]{4.6cm}{\center $f \in UC_{\text{loc}}(bv_1(\mathbb R),bv_p(\mathbb R))$\\ for $1 \leq p <+\infty$} & $\Leftarrow$ & \parbox[c]{12cm}{\center for every $r>0$ and $\varepsilon>0$ there exists $\delta>0$ such that for all $a,b,c,d \in [-r,r]$ with $a\neq b$, $c\neq d$ and $\abs{a-c}\leq \delta$, $\abs{b-d}\leq \delta$ we have
$\abs[\Bigg]{\dfrac{(f(a)-f(b))(a-b)}{\abs{a-b}^{1+\frac{1}{p}}} - \dfrac{(f(c)-f(d))(c-d)}{\abs{c-d}^{1+\frac{1}{p}}} } \leq \varepsilon$} \\[4mm] \midrule 
\parbox[c]{4.6cm}{\center $f \in UC_{\text{loc}}(bv_1(\mathbb R),bv_1(\mathbb R))$} & $\Leftrightarrow$ & \parbox[c]{10.5cm}{\center $f$ is continuously differentiable} \\[2mm] \midrule 
\parbox[c]{4.4cm}{\center $f \in UC_{\text{loc}}(bv_p(E),bv_q(E))$\\ for $1< p \leq q<+\infty$} & $\Leftarrow$ & \parbox[c]{10.5cm}{\center the map $u \mapsto f(u)-f(0)$ is continuous and $\mathbb R$-linear} \\[3.5mm] \midrule 
\parbox[c]{4.4cm}{\center $f \in UC_{\text{loc}}(bv_p(E),c_0(E))$\\ for $1 \leq p <+\infty$} & $\Leftrightarrow$ & $f$ is the zero map \\[3mm] \midrule 
\parbox[c]{4.4cm}{\center $f \in UC_{\text{loc}}(bv_p(E),l^q(E))$\\ for $1 \leq p,q <+\infty$} & $\Leftrightarrow$ & $f$ is the zero map \\[3mm] \midrule 
$f \in UC_{\text{loc}}(bv_1(\hat E),l^\infty(\hat E))$ & $\Leftrightarrow$ & $f$ is locally uniformly continuous on $\hat E$\\ \midrule
\parbox[c]{4.4cm}{\center $f \in UC_{\text{loc}}(bv_p(E),l^\infty(E))$\\ for $1<p<+\infty$} & $\Leftrightarrow$ & $f$ is constant \\[3mm] \midrule
$f \in UC_{\text{loc}}(bv_1(\hat E),c(\hat E))$ & $\Leftrightarrow$ & $f$ is locally uniformly continuous on $\hat E$ \\ \midrule
\parbox[c]{4.4cm}{\center $f \in UC_{\text{loc}}(bv_p(E),c(E))$\\ for $1<p<+\infty$} & $\Leftrightarrow$ & $f$ is constant \\[2mm]
\bottomrule
\end{tabular}
\vspace{1.5mm}
\caption{Local uniform continuity for the composition operator $C_f$.}
\end{center}
\end{table}
\end{small}

\clearpage

\begin{small}
\begin{table}[hbt!]
\begin{center}
\begin{tabular}{@{}ccc@{}}
\toprule
\parbox[c]{4.4cm}{\center $f \in UC(c_0(E),bv_p(E))$\\ for $1 \leq p < +\infty$} & $\Leftrightarrow$ & $f$ is constant\\[3mm] \midrule 
\parbox[c]{4.4cm}{\center $f \in UC(c(E),bv_p(E))$\\ for $1\leq p < +\infty$} & $\Leftrightarrow$ & $f$ is constant\\[3mm] \midrule 
\parbox[c]{4.4cm}{\center $f \in UC(l^\infty(E),bv_p(E))$\\ for $1\leq p < +\infty$} & $\Leftrightarrow$ & $f$ is constant\\[3mm] \midrule 
\parbox[c]{4.4cm}{\center $f \in UC(l^p(E),bv_q(E))$\\ for $1\leq p \leq q < +\infty$} & $\Leftrightarrow$ & \parbox[c]{11cm}{\center there exist $\delta>0$ and $L\geq 0$ such that\linebreak $\norm{f(u)-f(w)}\leq L\norm{u-w}^{\frac{p}{q}}$ for all $u,w \in E$ with $\norm{u-w}\leq \delta$} \\[3mm] \midrule
 \parbox[c]{4.4cm}{\center $f \in UC(l^p(E),bv_q(E))$\\ for $1\leq q < p < +\infty$} & $\Leftrightarrow$ & \parbox[c]{12cm}{\center $f$ is constant} \\[3mm] \midrule
\parbox[c]{4.4cm}{\center $f \in UC(bv_p(E),bv_q(E))$\\ for $1\leq q < p <+\infty$} & $\Leftrightarrow$ & $f$ is constant \\[3mm] \midrule 
\parbox[c]{4.4cm}{\center $f \in UC(bv_1(\hat E),bv_p(\hat E))$\\ for $1\leq p < +\infty$} & $\Leftrightarrow$ & \parbox[c]{10.5cm}{\center the map $u \mapsto f(u)-f(0)$ is continuous and $\mathbb R$-linear} \\[4mm] \midrule 
\parbox[c]{4.4cm}{\center $f \in UC(bv_p(E),bv_q(E))$\\ for $1< p \leq q<+\infty$} & $\Leftrightarrow$ & \parbox[c]{10.5cm}{\center the map $u \mapsto f(u)-f(0)$ is continuous and $\mathbb R$-linear} \\[3.5mm] \midrule 
\parbox[c]{4.4cm}{\center $f \in UC(bv_p(E),c_0(E))$\\ for $1 \leq p <+\infty$} & $\Leftrightarrow$ & $f$ is the zero map \\[3mm] \midrule 
\parbox[c]{4.4cm}{\center $f \in UC(bv_p(E),l^q(E))$\\ for $1 \leq p,q <+\infty$} & $\Leftrightarrow$ & $f$ is the zero map \\[3mm] \midrule 
$f \in UC(bv_1(\hat E),l^\infty(\hat E))$ & $\Leftrightarrow$ & $f$ is uniformly continuous on $\hat E$\\ \midrule
\parbox[c]{4.4cm}{\center $f \in UC(bv_p(E),l^\infty(E))$\\ for $1<p<+\infty$} & $\Leftrightarrow$ & $f$ is constant \\[3mm] \midrule
$f \in UC(bv_1(\hat E),c(\hat E))$ & $\Leftrightarrow$ & $f$ is uniformly continuous on $\hat E$ \\ \midrule
\parbox[c]{4.4cm}{\center $f \in UC(bv_p(E),c(E))$\\ for $1<p<+\infty$} & $\Leftrightarrow$ & $f$ is constant \\[2mm]
\bottomrule
\end{tabular}
\vspace{1.5mm}
\caption{Uniform continuity for the composition operator $C_f$.}
\end{center}
\end{table}
\end{small}

\clearpage

\section{H\"older continuity}

In this section, we continue our study of the continuity properties of composition operators, now focusing on Lipschitz and H\"older continuity. As before, we distinguish between global and local cases. Interestingly, these two notions behave quite differently in some spaces while being nearly indistinguishable in others.

The definitions of various types of H\"older continuity for nonlinear operators between different normed spaces are well-known and closely follow those presented in Section~\ref{sec:composition_operators}. Therefore, we will not repeat them. As before, we omit the phrase ``with exponent $\alpha$'' whenever the exponent is clear from the context. It is also important to emphasize that the composition operator $C_f$ and its generator $f$ are defined on different spaces, which necessitates interpreting their H\"older continuity in distinct frameworks. For instance, in Theorem~\ref{prop:general_lip} this distinction applies to the pair $X$ and $Y$ for $C_f$ and the space $E$ for $f$. Since the relevant context will always be evident, we omit these distinctions to keep the presentation concise.

Let us start with a general result

\begin{proposition}\label{prop:general_lip}
Let $p,q \in [1,+\infty)$. Moreover, let $X \in \{l^p(E),c_0(E),c(E),l^\infty(E),bv_p(E)\}$ and $Y \in \{l^q(E),c_0(E),c(E),l^\infty(E),bv_q(E)\}$. Assume also that the composition operator $C_f$ maps $X$ into $Y$. If $C_f$ is H\"older continuous \textup(respectively, H\"older continuous on bounded sets\textup) with exponent $\alpha \in (0,1]$, then the function $f$  possesses the same H\"older continuity property.
\end{proposition}

\begin{proof}
(a) \emph{H\"older continuity on bounded sets}. Assume that the composition operator $C_f \colon X \to Y$ satisfies the H\"older condition on bounded sets with exponent $\alpha \in (0,1]$. Fix $r>0$ and let $L_r\geq 0$ be the H\"older constant of $C_f$ restricted to the closed ball $\ball_X(0,r)$. For any $u,w \in \ball_E(0,r)$ define the sequences $x:=(u,0,0,0,\ldots)$, $y:=(w,0,0,0,\ldots)$ if $X\in \{l^p(E),c_0(E),c(E),l^\infty(E)\}$ and $x:=(u,u,u,u,\ldots)$, $y:=(w,w,w,w,\ldots)$ if $X=bv_p(E)$. Since $x,y \in \ball_X(0,r)$, applying the H\"older continuity of $C_f$ gives
\[
 \norm{f(u)-f(w)}\leq \norm{C_f(x)-C_f(y)}_Y \leq L_r\norm{x-y}_X^\alpha = L_r\norm{u-w}^\alpha.
\]
Thus, $f$ is H\"older continuous on bounded sets with exponent $\alpha$.

(b) \emph{H\"older continuity}. The argument proceeds in the same manner as in part~(a). Hence, it is omitted.
\end{proof}

\subsection{H\"older continuity on bounded sets}
\label{sec:local_holder}

In certain cases, the converse of the above proposition holds, particularly when the target space of the composition operator is $c(E)$ or $l^\infty(E)$.

\begin{theorem}\label{thm:local_holder_c_l_infty_bv}
Let $E$ be a Banach space and let $X \in \{c(E), l^\infty(E)\}$. Moreover, assume that the composition operator $C_f$ maps $bv_1(E)$ into $X$. Then, $C_f$ is H\"older continuous on bounded sets with exponent $\alpha \in (0,1]$ if and only if $f$ is.
\end{theorem}

\begin{proof}
By Proposition~\ref{prop:general_lip}, it is enough to prove the sufficiency part. Let $r>0$ and $L_r\geq 0$ be such that $\norm{f(u)-f(w)} \leq L_r \norm{u-w}^\alpha$ for all $u,w \in \ball_{E}(0,r)$. If $x,y\in \ball_{bv_1}(0,r)$, then for each $n \in \mathbb N$ we have $x(n),y(n) \in  \ball_{E}(0,r)$. Thus, we  obtain
\begin{align*}
\norm{C_f(x)-C_f(y)}_\infty &= \sup_{n \in \mathbb N} \norm{f(x(n))-f(y(n))} \leq \sup_{n \in \mathbb N} L_r\norm{x(n)-y(n)}^\alpha\\
 & = L_r \norm{x-y}_{\infty}^\alpha \leq L_r \norm{x-y}^\alpha_{bv_1}.
\end{align*}
This concludes the proof.
\end{proof}

If the composition operator acts from the space $l^p(E)$ to $bv_q(E)$, the situation becomes more complicated and depends on the relationship between the exponent $\alpha$ and the parameters $p$ and $q$. 

\begin{theorem}\label{thm:holder_of_lp_bvq_part1}
Let $1 \leq p \leq q < +\infty$ and $\frac{p}{q}\leq \alpha \leq 1$. Moreover, assume that the composition operator $C_f$ maps $l^p(E)$ into $bv_q(E)$. Then, $C_f$ is H\"older continuous on bounded sets with exponent $\alpha$ if and only if $f$ is.
\end{theorem}

\begin{proof}
In view of Proposition~\ref{prop:general_lip}, it suffices to prove the ``if'' part. Fix $r>0$ and consider any $x,y \in \ball_{l^p}(0,r)$. Let $L_r \geq 0$ denote the H\"older constant of $f$ restricted to $\ball_E(0,r)$. Since $\norm{x(n)}\leq r$ and $\norm{y(n)}\leq r$ for all $n \in \mathbb N$, we have
\[
 \norm{f(x(1))-f(y(1))}\leq L_r \norm{x(1)-y(1)}^\alpha \leq L_r\norm{x-y}_{l^p}^\alpha.
\]
Furthermore,
\begin{align*}
 &\Biggl(\sum_{n=1}^\infty \norm{f(x(n+1))-f(y(n+1))-f(x(n))+f(y(n))}^q \Biggr)^{\frac{1}{q}}\\
 &\quad \leq \Biggl(2^{q+1}\sum_{n=1}^\infty \norm{f(x(n))-f(y(n))}^q \Biggr)^{\frac{1}{q}} \leq 2^{1+\frac{1}{q}}L_r \Biggl(\sum_{n=1}^\infty \norm{x(n)-y(n)}^{\alpha q}\Biggr)^{\frac{1}{q}}.
\end{align*}
Since $p \leq \alpha q$, we have
\[
 \Biggl(\sum_{n=1}^\infty \norm{x(n)-y(n)}^{\alpha q}\Biggr)^{\frac{1}{q}} = \Biggl(\sum_{n=1}^\infty \norm{x(n)-y(n)}^{\alpha q}\Biggr)^{\frac{\alpha}{\alpha q}} \leq \Biggl(\sum_{n=1}^\infty \norm{x(n)-y(n)}^{p}\Biggr)^{\frac{\alpha}{p}}
\]
(cf.~\cite{HLP}*{Theorem~19, p.~28}). Thus, we obtain
\[
 \Biggl(\sum_{n=1}^\infty \norm{f(x(n+1))-f(y(n+1))-f(x(n))+f(y(n))}^q \Biggr)^{\frac{1}{q}} \leq 2^{1+\frac{1}{q}}L_r \norm{x-y}_{l^p}^\alpha.
\]
Finally, we conclude that $\norm{C_f(x)-C_f(y)}_{bv_q} \leq 5L_r \norm{x-y}_{l^p}^\alpha$. This completes the proof. 
\end{proof}

If $0< \alpha < \frac{p}{q}\leq 1$ we will provide a sufficient condition only, and we do not know whether it is also a necessary one. 

\begin{theorem}\label{thm:holder_of_lp_bvq_part2}
Let $1 \leq p \leq q < +\infty$ and let $0< \alpha \leq \frac{p}{q}\leq 1$. If $f$ is H\"older continuous on bounded sets with exponent $\alpha$ and there exists a closed ball $\ball_E(0,R)$, where $R>0$, on which it satisfies the H\"older condition with exponent $\frac{p}{q}$, then the composition operator $C_f$ maps $l^p(E)$ into $bv_q(E)$ and is H\"older continuous on bounded sets with exponent $\alpha$.
\end{theorem}

\begin{proof}
The fact that the composition operator $C_f$ maps $l^p(E)$ into $bv_q(E)$ follows directly from~\cite{BK}*{Theorem~17}.

Let $r>0$ be fixed. Without loss of generality, we may assume that $r> R$. Denote by $L_r\geq 0$ the H\"older constant of $f$ corresponding to the ball $\ball_E(0,r)$. Similarly, let $\Lambda_R \geq 0$ represent the  H\"older constant of $f$ restricted to the ball $\ball_E(0,R)$, where the exponent is $\frac{p}{q}$. Note that $L_r$ corresponds to the exponent $\alpha$. Now consider any $x,y \in \ball_{l^p}(0,r)$. We partition the set of indices $\mathbb N$ into four disjoint sets based on the norms of $x(n)$ and $y(n)$:
\begin{align*}
   I_1&:=\dset[\big]{n \in \mathbb N}{\text{$\norm{x(n)}\geq R$ and $\norm{y(n)}\geq R$}},\\[1mm]
	I_2&:=\dset[\big]{n \in \mathbb N}{\text{$\norm{x(n)}< R$ and $\norm{y(n)}\geq R$}},\\[1mm]
	I_3&:=\dset[\big]{n \in \mathbb N}{\text{$\norm{x(n)}\geq R$ and $\norm{y(n)} < R$}},\\[1mm]
	I_4&:=\dset[\big]{n \in \mathbb N}{\text{$\norm{x(n)} < R$ and $\norm{y(n)}< R$}}.
\end{align*} 
The sets $I_1,I_2,I_3$ are clearly finite (possibly empty), while $I_4$ contains infinitely many positive integers. Let $\abs{J}$ denote the cardinality of $J$, provided it is finite. Then, for any $j\in \{1,2,3\}$ we have
\[
 r^p \geq \sum_{n=1}^\infty \norm{z(n)}^p \geq \sum_{n \in I_j} \norm{z(n)}^p \geq \abs{I_j}\cdot R^p,
\]
where $z$ represents $x$ if $j=1,3$, and $y$ if $j=2$. Consequently, we obtain the upper bound for the cardinalities of the sets $I_1,I_2$ and $I_3$:
\[
\abs{I_j} \leq \frac{r^p}{R^p} \ \text{for $j=1,2,3$}.
\]
Notably, this bound is independent of the specific choice of points within the ball $\ball_{l^p}(0,r)$.

Let us now proceed with estimating the sum
\[
 \Biggl(\sum_{n=1}^\infty \norm{f(x(n+1))-f(y(n+1))-f(x(n))+f(y(n))}^q \Biggr)^{\frac{1}{q}}.
\]
We have
\begin{align*}
 &\Biggl(\sum_{n=1}^\infty \norm{f(x(n+1))-f(y(n+1))-f(x(n))+f(y(n))}^q \Biggr)^{\frac{1}{q}}\\
 &\qquad \leq \Biggl(2^{q+1}\sum_{n=1}^\infty \norm{f(x(n))-f(y(n))}^q \Biggr)^{\frac{1}{q}}\\
 &\qquad = \Biggl(2^{q+1}\sum_{j=1}^4\sum_{n \in I_j} \norm{f(x(n))-f(y(n))}^q \Biggr)^{\frac{1}{q}}\\
 &\qquad \leq 4\sum_{j=1}^4 \Biggl(\sum_{n \in I_j} \norm{f(x(n))-f(y(n))}^q \Biggr)^{\frac{1}{q}}.
\end{align*}
If any of the sets $I_j$ is empty, we define the corresponding sum as zero by convention.

We can now focus on each sum over the set $I_j$ separately. We begin with the sum corresponding to $I_1$. First, by applying the H\"older continuity of $f$ with exponent $\alpha$, then using Jensen's inequality, and finally utilizing the estimate for the cardinality of $I_1$, we obtain
\begin{align*}
 \Biggl(\sum_{n \in I_1} \norm{f(x(n))-f(y(n))}^q \Biggr)^{\frac{1}{q}} & \leq  \Biggl(L_r^q \sum_{n \in I_1} \norm{x(n)-y(n)}^{\alpha q} \Biggr)^{\frac{p}{q\alpha} \cdot \frac{\alpha}{p}}\\
& \leq L_r \Biggl(\abs{I_1}^{\frac{p}{\alpha q}-1}\cdot \sum_{n \in I_1} \norm{x(n)-y(n)}^{p} \Biggr)^{\frac{\alpha}{p}}\\
& \leq L_r \cdot \biggl(\frac{r}{R}\biggr)^{\frac{p}{q} - \alpha}\cdot \Biggl(\sum_{n \in I_1} \norm{x(n)-y(n)}^{p} \Biggr)^{ \frac{\alpha}{p}}\\
& \leq L_r\cdot \biggl(\frac{r}{R}\biggr)^{\frac{p}{q} - \alpha}\cdot \norm{x-y}_{l^p}^\alpha.
\end{align*} 

We now proceed to estimate the second sum. Fix $n \in I_2$, and for $t \in [0,1]$ define $w_t(n):=(1-t)x(n)+ty(n)$. Since $\norm{w_0(n)}=\norm{x(n)}<R$ and $\norm{w_1(n)}=\norm{y(n)}\geq R$, there exists a value $t_\ast:=t_\ast(n) \in (0,1]$ such that $\norm{w_{t_\ast}(n)}=R$. For simplicity, let us set $z(n):=w_{t_\ast}(n)$ for all $n \in I_2$. Then, by applying the Minkowski inequality, we obtain
\begin{align*}
 &\Biggl(\sum_{n \in I_2} \norm{f(x(n))-f(y(n))}^q \Biggr)^{\frac{1}{q}}\\
 &\qquad\leq  \Biggl(\sum_{n \in I_2} \norm{f(x(n))-f(z(n))}^q \Biggr)^{\frac{1}{q}} +  \Biggl(\sum_{n \in I_2} \norm{f(z(n))-f(y(n))}^q \Biggr)^{\frac{1}{q}}.
\end{align*} 
Next, by employing a similar reasoning to that used for the set $I_1$, we can estimate the second term appearing on the right-hand side of the inequality above. We have
\[
 \Biggl(\sum_{n \in I_2} \norm{f(z(n))-f(y(n))}^q \Biggr)^{\frac{1}{q}} \leq L_r \cdot \biggl(\frac{r}{R}\biggr)^{\frac{p}{q} - \alpha}\cdot \Biggl(\sum_{n \in I_2} \norm{z(n)-y(n)}^{p} \Biggr)^{ \frac{\alpha}{p}}.
\]
Since $\norm{z(n)-y(n)}=(1-t_{\ast})\norm{x(n)-y(n)}$, we obtain
\[
  \Biggl(\sum_{n \in I_2} \norm{f(z(n))-f(y(n))}^q \Biggr)^{\frac{1}{q}} \leq  L_r \cdot \biggl(\frac{r}{R}\biggr)^{\frac{p}{q} - \alpha}\cdot \Biggl(\sum_{n \in I_2} \norm{x(n)-y(n)}^{p} \Biggr)^{ \frac{\alpha}{p}}.
\]
As for the first term, applying the H\"older continuity of $f$ restricted to the ball $\ball_E(0,R)$ with exponent $\frac{p}{q}$, we arrive at the estimate:
\begin{align*}
& \Biggl(\sum_{n \in I_2} \norm{f(x(n))-f(z(n))}^q \Biggr)^{\frac{1}{q}}\\
&\qquad  \leq \Lambda_R \Biggl(\sum_{n \in I_2} \norm{x(n)-z(n)}^{p} \Biggr)^{\frac{1}{q}} \leq \Lambda_R \Biggl(\sum_{n \in I_2} \norm{x(n)-y(n)}^{p} \Biggr)^{\frac{1}{q}}\\
&\qquad \leq \Lambda_R \norm{x-y}_{l^p}^{\frac{p}{q}}=\Lambda_R \norm{x-y}_{l^p}^{\alpha} \cdot \norm{x-y}_{l^p}^{\frac{p}{q}-\alpha} \leq \Lambda_R \cdot (2r)^{\frac{p}{q}-\alpha} \cdot \norm{x-y}_{l^p}^{\alpha}.
\end{align*}
Combining both terms, we conclude that
\[
\Biggl(\sum_{n \in I_2} \norm{f(x(n))-f(y(n))}^q \Biggr)^{\frac{1}{q}} \leq \biggl( L_r \cdot \biggl(\frac{r}{R}\biggr)^{\frac{p}{q} - \alpha}+ \Lambda_R \cdot (2r)^{\frac{p}{q}-\alpha}\biggr)\ \norm{x-y}_{l^p}^\alpha.
\]

We will omit the details of the estimates for the last two sums, as they can be derived using a similar approach to the one previously applied. Specifically, we obtain
\[
  \Biggl(\sum_{n \in I_3} \norm{f(x(n))-f(y(n))}^q \Biggr)^{\frac{1}{q}} \leq \biggl( L_r \cdot \biggl(\frac{r}{R}\biggr)^{\frac{p}{q} - \alpha}+ \Lambda_R \cdot (2r)^{\frac{p}{q}-\alpha}\biggr)\ \norm{x-y}_{l^p}^\alpha
\]
and
\[	
	\Biggl(\sum_{n \in I_4} \norm{f(x(n))-f(y(n))}^q \Biggr)^{\frac{1}{q}} \leq \Lambda_R \cdot (2r)^{\frac{p}{q}-\alpha} \cdot \norm{x-y}_{l^p}^\alpha.
\]
(Note that, in establishing the upper bound for the first term in the sum corresponding to the set $I_2$, we could have alternatively applied Jensen's inequality along with the H\"older continuity of $f$ with exponent $\alpha$, rather than relying on the H\"older continuity of $f$ with exponent $\frac{p}{q}$ on the ball $\ball_E(0,R)$. However, in this instance, a straightforward application of Jensen’s inequality may not suffice to derive the desired estimate, as the set $I_4$ is infinite. To ensure the argument holds we must assume the stronger H\"older continuity property for $f$.)  

By summing all the estimates and considering the fact that
\[
\norm{f(x(1))-f(y(1))}\leq L_r\norm{x(1)-y(1)}^\alpha \leq L_r\norm{x-y}_{l^p}^\alpha,
\]
we ultimately obtain the desired inequality $\norm{C_f(x)-C_f(y)}_{bv_q} \leq \Gamma_r \norm{x-y}_{l^p}^\alpha$, where
\[
\Gamma_r:=L_r + 12L_r \cdot \biggl(\frac{r}{R}\biggr)^{\frac{p}{q} - \alpha} + 12 \Lambda_R \cdot (2r)^{\frac{p}{q}-\alpha}.
\]
Thus, we conclude that the composition operator $C_f \colon l^p(E) \to bv_q(E)$ is H\"older continuous on bounded sets with exponent $\alpha$.
\end{proof}

As previously noted, we are unsure whether the conditions outlined in Theorem~\ref{thm:holder_of_lp_bvq_part2} for the generator $f$ are also necessary for the composition operator $C_f \colon l^p(E) \to bv_q(E)$ to be H\"older continuous with exponent $\alpha$, where $0< \alpha < \frac{p}{q}\leq 1$. It is clear that merely assuming that $f \in \Lip_{\text{bnd}}^\alpha(E)$ is insufficient, as in this case the composition operator may not even be well-defined as a map between $l^p(E)$ and $bv_q(E)$  (cf.~\cite{BK}*{Theorem~17}); to see this take $E:=\mathbb R$ and consider the power function $u \mapsto \abs{u}^\alpha$. Moreover, we do not know whether a result analogous to Theorems~\ref{thm:holder_of_lp_bvq_part1} and~\ref{thm:holder_of_lp_bvq_part2} can be proved when $\alpha \leq 1<\frac{p}{q}$. However, as the following example demonstrates, non-trivial Lipschitz continuous composition operators exist in this case.  

\begin{example}
Let $q \in [1,+\infty)$ and set $p:=2q$. Moreover, assume that $E$ is a commutative normed algebra. Define the map $f \colon E \to E$ by $f(u)=u^2$. To demonstrate that $f$ generates a composition operator $C_f \colon l^p(E) \to bv_q(E)$ that is Lipschitz continuous on bounded sets, we fix an arbitrary $r>0$ and consider two points $x,y \in \ball_{l^p}(0,r)$. First observe that
\[
 \norm{f(x(1))-f(y(1))} = \norm[\big]{x^2(1)-y^2(1)} \leq \norm{x(1)-y(1)}\cdot \norm{x(1)+y(1)} \leq 2r\norm{x-y}_{l^p}.
\]
Further, using the Cauchy--Schwarz inequality, we obtain
\begin{align*}
 \sum_{n=1}^\infty \norm[\big]{x^2(n)-y^2(n)}^q & \leq \sum_{n=1}^\infty \norm{x(n)-y(n)}^q\cdot \norm{x(n)+y(n)}^q\\
 & \leq \Biggl(\sum_{n=1}^\infty \norm{x(n)-y(n)}^{2q}\Biggr)^{\frac{1}{2}} \cdot \Biggl(\sum_{n=1}^\infty \norm{x(n)+y(n)}^{2q}\Biggr)^{\frac{1}{2}}\\
& =  \Biggl(\sum_{n=1}^\infty \norm{x(n)-y(n)}^{p}\Biggr)^{\frac{q}{p}} \cdot \Biggl(\sum_{n=1}^\infty \norm{x(n)+y(n)}^{p}\Biggr)^{\frac{q}{p}}.
\end{align*} 
Thus, we conclude that
\[
\Biggl(\sum_{n=1}^\infty \norm[\big]{x^2(n)-y^2(n)}^q\Biggr)^{\frac{1}{q}}\leq \norm{x-y}_{l^p} \cdot \norm{x+y}_{l^p} \leq 2r\norm{x-y}_{l^p}.
\]
By the triangle inequality and the previous bounds, we finally get
\begin{align*}
& \Biggl(\sum_{n=1}^\infty \norm{f(x(n+1))-f(y(n+1))-f(x(n))+f(y(n))}^q \Biggr)^{\frac{1}{q}}\\
& \hspace{2cm}\leq 2^{1+\frac{1}{q}}  \Biggl(\sum_{n=1}^\infty \norm{f(x(n))-f(y(n))}^q \Biggr)^{\frac{1}{q}} \leq 8r \norm{x-y}_{l^p}.
\end{align*} 
This proves not only that the composition operator $C_f$ maps $l^p(E)$ into $bv_q(E)$, but also that $\norm{C_f(x)-C_f(y)}_{bv_q}\leq 10r\norm{x-y}_{l^p}$ for $x,y \in  \ball_{l^p}(0,r)$. 
\end{example}

Finally, we focus on the local Lipschitz continuity of composition operators acting between $bv_1(E)$ and $bv_p(E)$. We will provide only sufficient conditions. With the exception of the case when $p=1$ and $E:=\mathbb R$, we do not know whether these conditions are also necessary. 

\begin{theorem}\label{thm:LH1_bv_1_bv_p}
Let $p \in [1,+\infty)$ and let $E$ be a Banach space. If the generator $f \colon E \to E$ is  Fr\'echet differentiable with the derivative $f' \colon E \to \mathcal{L}(E)$ being H\"older continuous on bounded sets with exponent $\frac{1}{p}$, then the composition operator $C_f$ maps $bv_1(E)$ into $bv_p(E)$ and is Lipschitz continuous on bounded sets.   
\end{theorem}

\begin{proof}
First, observe that since the derivative $f'$ is H\"older continuous on bounded sets, it is also locally bounded. Consequently, $f$ satisfies a Lipschitz condition on bounded subsets of $E$  (cf. the proof of Theorem~\ref{thm:LUC_bv_1_bv_p}). In particular, by~\cite{BK}*{Theorem~20}, the composition operator $C_f$ maps $bv_1(E)$ into $bv_p(E)$.   

In the proof of Lipschitz continuity on bounded sets of $C_f$, we will again use the Lagrange-type formula~\eqref{eq:mean_value} and loosely follow the argument presented in~\cite{DN}*{Theorem~6.68}. Let $r>0$ be fixed, and let $L_r \geq 0$ denote the H\"older constant of the derivative $f'$ restricted to the closed ball $\ball_E(0,r)$. Additionally, set $M_r:=\sup_{\norm{u}\leq r}\norm{f'(u)}_{\mathcal L}$. Since $f'$ is H\"older continuous on bounded sets, the constant $M_r$ is finite. Now, if $x,y \in \ball_{bv_1}(0,r)$, then for each $n \in \mathbb N$ both $x(n)$ and $y(n)$ lie in $\ball_E(0,r)$. Therefore, we can estimate the difference of the first terms as follows:
\begin{align*}
 \norm{f(x(1))-f(y(1))}&=\norm[\bigg]{\int_0^1 f'\bigl(y(1)+t(x(1)-y(1))\bigr)(x(1)-y(1)) \textup dt}\\
&\leq \int_0^1 \norm[\big]{f'\bigl(y(1)+t(x(1)-y(1))\bigr)}_{\mathcal L} \textup dt \cdot \norm{x(1)-y(1)} \leq M_r \norm{x(1)-y(1)}.
\end{align*}
Furthermore, for any $n \in \mathbb N$ we have
\begin{align*}
 &\norm[\Big]{\bigl[f(x(n+1))-f(y(n+1))\bigr] - \bigl[f(x(n))-f(y(n))\bigr]}^p\\
 &\quad = \norm[\bigg]{\int_0^1 f'\bigl(y(n+1)+t(x(n+1)-y(n+1))\bigr)(x(n+1)-y(n+1))\textup dt\\
 &\quad \qquad - \int_0^1 f'\bigl(y(n)+t(x(n)-y(n))\bigr)(x(n)-y(n))\textup dt}^p\\
 &\quad = \norm[\bigg]{\int_0^1 \bigl[f'\bigl(y(n+1)+t(x(n+1)-y(n+1))\bigr)\\
 &\hspace{4cm} - f'\bigl(y(n)+t(x(n)-y(n))\bigr)\bigr](x(n+1)-y(n+1))\textup dt\\
 &\quad \qquad + \int_0^1 f'\bigl(y(n)+t(x(n)-y(n))\bigr)\bigl[(x(n+1)-x(n))-(y(n+1)-y(n))\bigr]\textup dt}^p.
\end{align*}
Applying the triangle inequality and then Jensen's inequality, we thus obtain
\begin{align*}
 &\norm[\Big]{\bigl[f(x(n+1))-f(y(n+1))\bigr] - \bigl[f(x(n))-f(y(n))\bigr]}^p\\
 &\qquad \leq 2^p I_1(n) \cdot \norm{x(n+1)-y(n+1)}^p + 2^p I_2(n) \cdot \norm[\big]{(x(n+1)-x(n))-(y(n+1)-y(n))}^p,
\end{align*}
where
\begin{align*}
 I_1(n)&:= \int_0^1 \norm[\Big]{f'\bigl(y(n+1)+t(x(n+1)-y(n+1))\bigr) - f'\bigl(y(n)+t(x(n)-y(n))\bigr)}_{\mathcal L}^p\textup dt  
\end{align*}
and
\begin{align*}
 I_2(n):=\int_0^1 \norm[\Big]{f'\bigl(y(n)+t(x(n)-y(n))\bigr)}_{\mathcal L}^p \textup dt.
\end{align*}

Regarding the integral $I_2(n)$, we can simply bound it by $M_r^p$, that is,  $I_2(n)\leq M_r^p$ for all $n \in \mathbb N$. On the other hand, to estimate $I_1(n)$ we will use the H\"older continuity of $f'$. For a fixed $n \in \mathbb N$ and each $t \in [0,1]$ we have
\begin{align*}
&\norm[\Big]{f'\bigl(y(n+1)+t(x(n+1)-y(n+1))\bigr)- f'\bigl(y(n)+t(x(n)-y(n))\bigr)}_{\mathcal L}^p\\
&\qquad \leq L_r^p \norm[\big]{y(n+1)-y(n)+t(x(n+1)-x(n)-y(n+1)+y(n))}\\[1mm]
&\qquad \leq L_r^p\norm[\big]{y(n+1)-y(n)} + L_r^p \norm[\big]{x(n+1)-y(n+1)-x(n)+y(n)}.
\end{align*}  
Thus, we find that
\[
 I_1(n)\leq L_r^p\norm[\big]{y(n+1)-y(n)} + L_r^p \norm[\big]{x(n+1)-y(n+1)-x(n)+y(n)}.
\]
Summing over $n$, we obtain
\[
 \sum_{n=1}^\infty I_1(n) \leq L_r^p \norm{y}_{bv_1} + L_r^p \norm{x-y}_{bv_1} \leq 3r L_r^p.
\]
Consequently, we have
\begin{align*}
&\sum_{n=1}^\infty \norm[\Big]{\bigl[f(x(n+1))-f(y(n+1))\bigr] - \bigl[f(x(n))-f(y(n))\bigr]}^p\\
&\hspace{6cm} \leq 3\cdot 2^p r L_r^p \norm{x-y}_{\infty}^p + 2^p M_r^p \norm{x-y}_{bv_p}^p.
\end{align*}
Therefore,
\begin{align*}
&\Biggl(\sum_{n=1}^\infty \norm[\Big]{\bigl[f(x(n+1))-f(y(n+1))\bigr] - \bigl[f(x(n))-f(y(n))\bigr]}^p\Biggr)^{\frac{1}{p}}\\
&\qquad \leq  2\cdot 3^{\frac{1}{p}} r^{\frac{1}{p}} L_r \norm{x-y}_{\infty} + 2 M_r \norm{x-y}_{bv_p} \leq 6r^{\frac{1}{p}} L_r \norm{x-y}_{bv_1} + 2M_r \norm{x-y}_{bv_1}.
\end{align*}
Finally, combining all estimates, we obtain the desired inequality $\norm{C_f(x)-C_f(y)}_{bv_p} \leq \Lambda_r \norm{x-y}_{bv_1}$, where $\Lambda_r:= 6r^{\frac{1}{p}} L_r + 3M_r$. The proof is now complete.
\end{proof}

In the special case when $E:=\mathbb R$ and $p=1$ Theorem~\ref{thm:LH1_bv_1_bv_p} can be inverted. 

\begin{theorem}\label{thm:LH1_bv_1_bv_p_ver2}
A composition operator $C_f$ maps $bv_1(\mathbb R)$ into itself and is Lipschitz continuous on bounded sets if and only if the derivative $f'$ of its generator $f$ satisfies a Lipschitz condition on bounded sets.
\end{theorem}

\begin{proof}
In view of Theorem~\ref{thm:LH1_bv_1_bv_p}, it suffices to prove only the necessity part. To this end, we will employ an approach inspired by~\cite{reinwand}*{Theorem~5.1.21} -- cf. also the proof of Theorem~\ref{thm:LUC_C1}.

Assume that the composition operator $C_f$ maps the space $bv_1(\mathbb R)$  into itself and is Lipschitz continuous on bounded sets. From the characterization of acting conditions for $C_f \colon bv_1(\mathbb R) \to bv_1(\mathbb R)$ it follows that $f$ is Lipschitz continuous on every compact subset of $\mathbb R$ (see~\cite{BK}*{Theorem~20}). Consequently, $f$ is absolutely continuous on every non-empty and bounded interval. Let $D_f$ denote the set of points in $\mathbb R$ where $f$ is differentiable.

Fix $r\geq \frac{1}{16}$ and consider two points $u,w \in [-r,r]\cap D_f$. For a number $\rho \in (0,r]$ choose a positive integer $m$ such that $m \leq \frac{16r}{\rho}\leq 2m$ and define the sequences $x_\rho,y_\rho \colon \mathbb N \to \mathbb R$ as in the proof of Theorem~\ref{thm:LUC_C1}, that is, 
\begin{align*}
 x_\rho:&=(u,u+\rho,u,u+\rho,\ldots,u,u+\rho,u,u,u,\ldots),\\
 y_\rho:&=(w,w+\rho,w,w+\rho,\ldots,w,w+\rho,w,w,w,\ldots),
\end{align*}
where the last term of the form $u+\rho$ or $w+\rho$ appears at the $2m$-th position. Then, $x_\rho,y_\rho \in \ball_{bv_1}(0,33r)$ and
\begin{align*}
&2m\abs[\Big]{[f(u+\rho)-f(u)] - [f(w+\rho)-f(w)]}\\
&\qquad \leq \norm{C_f(x_\rho)-C_f(y_\rho)}_{bv_1} \leq L_r\norm{x_\rho-y_\rho}_{bv_1}=L_r\abs{u-w},
\end{align*}
where $L_r \geq 0$ is a Lipschitz constant corresponding to the restriction of $C_f$ to the closed ball $\ball_{bv_1}(0,33r)$. Consequently, for each $\rho \in (0,r]$, we obtain
\[
 \abs[\bigg]{\frac{f(u+\rho)-f(u)}{\rho} - \frac{f(w+\rho)-f(w)}{\rho}}\leq L_r\abs{u-w}.
\]
Taking the limit as $\rho \to 0^+$, we get $\abs{f'(u)-f'(w)} \leq L_r\abs{u-w}$. This implies that the restriction of $f'$ to the set $[-r,r]\cap D_f$ is Lipschitz continuous. To conclude it suffices to apply the following result: if the derivative $g'$ of an absolutely continuous function $g \colon [a,b] \to \mathbb R$ is Lipschitz continuous on the set where it exists, then it is Lipschitz continuous on the entire interval $[a,b]$ (see~\cite{reinwand}*{Lemma~1.1.27~(b)}). 
\end{proof}

\subsection{H\"older continuity} 

We initiate our study of a H\"older condition for composition operators with two results that closely resemble Theorems~\ref{thm:local_holder_c_l_infty_bv} and~\ref{thm:holder_of_lp_bvq_part1} from the previous subsection. For brevity, we state them without proof.

\begin{theorem}
Let $E$ be a Banach space and let $X \in \{c(E), l^\infty(E)\}$. Moreover, assume that the composition operator $C_f$ maps $bv_1(E)$ into $X$. Then, $C_f$ is H\"older continuous with exponent $\alpha \in (0,1]$ if and only if $f$ is.
\end{theorem}

\begin{theorem}
Let $1 \leq p \leq q < +\infty$ and let $\frac{p}{q}\leq \alpha \leq 1$. Moreover, assume that the composition operator $C_f$ maps $l^p(E)$ into $bv_q(E)$. Then, $C_f$ is H\"older continuous with exponent $\alpha$ if and only if $f$ is.
\end{theorem}

In Theorem~\ref{thm:holder_of_lp_bvq_part2} and the subsequent example, we saw that when the exponent $\alpha$ is smaller than $\frac{p}{q}$, there exist non-constant composition operators between $l^p(E)$ and $bv_q(E)$ that are H\"older continuous on bounded sets. However, the situation changes drastically when we deal with the (global) H\"older condition. 

\begin{theorem}
Let $p,q \in [1,+\infty)$ and let $\alpha \in (0,1]\cap (0,\frac{p}{q})$. Moreover, assume that the composition operator $C_f$ maps $l^p(E)$ into $bv_q(E)$. Then, $C_f$ is H\"older continuous with exponent $\alpha$ if and only if the map $f$ is constant.
\end{theorem}

\begin{proof}
Clearly, we need to prove the ``only if'' part. Assume that the composition operator $C_f$ is H\"older continuous with the constant $L \geq 0$ and exponent $\alpha$. Consider an arbitrary point $u \in E$ and for each $n \in \mathbb N$ set $x_n:=(0,u,0,u,0\ldots,0,u,0,0,0,\ldots)$, where the last occurrence of $u$ appears at the $2n$-th position. Clearly, $x_n \in l^p(E)$ for all $n \in \mathbb N$. Moreover, we have
\begin{align*}
 n^{\frac{1}{q}}\norm{f(u)-f(0)} \leq \norm{C_f(x_n)-C_f(0)}_{bv_q} \leq L\norm{x_n}^\alpha_{l^p} = Ln^{\frac{\alpha}{p}}\norm{u}^\alpha.
\end{align*}
Rearranging, this gives
\[
\norm{f(u)-f(0)}\leq Ln^{\frac{\alpha}{p}-\frac{1}{q}}\norm{u}^\alpha.
\]
Taking the limit in the above estimate as $n \to +\infty$, we obtain $\norm{f(u)-f(0)}=0$, because $\alpha<\frac{p}{q}$. This shows that $f(u)=f(0)$ for every $u \in E$ and completes the proof.
\end{proof}

We conclude our investigation of H\"older continuous composition operators with a result concerning $bv_p$-spaces that refines Theorem~\ref{thm:UC_bvp_bvq}.

\begin{theorem}
Let $1 \leq p \leq q < +\infty$. Moreover, assume that $E$ is a Banach space when $p=1$, and a normed space when $p>1$. Assume also that the composition operator $C_f$ maps $bv_p(E)$ into $bv_q(E)$. Then, $C_f$ is H\"older continuous with exponent $\alpha \in (0,1)$ if and only if $f$ is a constant map.
\end{theorem}

\begin{proof}
Clearly, only the necessity part requires proof. Assume that the composition operator $C_f$ is H\"older continuous with exponent $\alpha \in (0,1)$. This immediately implies that $C_f$ is uniformly continuous. Consequently, by Theorem~\ref{thm:UC_bvp_bvq}, the map $\varphi(u):=f(u)-f(0)$ is an $\mathbb R$-linear endomorphism of $E$. On the other hand, Proposition~\ref{prop:general_lip} ensures that $f$ satisfies the H\"older condition with exponent $\alpha$. Thus, there exists a constant $L\geq 0$ such that $\norm{\varphi(u)}\leq L\norm{u}^\alpha$ for all $u \in E$. Now, let us choose a point $w \in E$ with $\norm{w}=1$. Then, for every $t>0$ we obtain $t\cdot \norm{\varphi(w)}=\norm{\varphi(tw)}\leq L\norm{tw}^\alpha=Lt^\alpha$. Rearranging, this gives $\norm{\varphi(w)}\leq Lt^{\alpha-1}$. Taking the limit in the above inequality as $t \to +\infty$, we conclude that $\varphi(w)=0$. Finally, using the homogeneity of $\varphi$, it follows that $\varphi(u)=0$ for all $u \in E$, which implies that $f$ is constant.
\end{proof}

\subsection{H\"older continuity: summary}
Let $E$ be a normed space over the real or complex field. If $E$ is a Banach space, we denote it by $\hat E$ to emphasise its completeness. For two sequence spaces $X$ and $Y$, we define the following classes of maps $f \colon E \to E$ that generate composition operators $C_f \colon X \to Y$ with certain H\"older continuity properties: 
\begin{align*}
 H_{\text{bnd}}^\alpha(X,Y)&:=\dset[\Big]{f\colon E \to E}{\parbox{12cm}{$C_f$ maps $X$ into $Y$ and is H\"older continuous on bounded sets with exponent $\alpha$}},\\
 H^\alpha(X,Y)&:=\dset{f\colon E \to E}{\text{$C_f$ maps $X$ into $Y$ and is H\"older continuous with exponent $\alpha$}}.
\end{align*} 
To fully characterize these classes, we will also incorporate results on continuity from the previous section, as well as results on acting conditions and boundedness from~\cite{BK}*{Sections~3 and~4}. Additionally, we clarify that if no constraints are specified on the exponent $\alpha$ below, we mean that the property holds for any $\alpha \in (0,1]$.

\clearpage

\begin{small}
\begin{table}[hbt!]
\begin{center}
\begin{tabular}{@{}ccc@{}}
\toprule
\parbox[c]{4.4cm}{\center $f \in H_{\text{bnd}}^\alpha(c_0(E),bv_p(E))$\\ for $1 \leq p < +\infty$} & $\Leftrightarrow$ & $f$ is constant\\[3mm] \midrule 
\parbox[c]{4.4cm}{\center $f \in H_{\text{bnd}}^\alpha(c(E),bv_p(E))$\\ for $1\leq p < +\infty$} & $\Leftrightarrow$ & $f$ is constant\\[3mm] \midrule 
\parbox[c]{4.4cm}{\center $f \in H_{\text{bnd}}^\alpha(l^\infty(E),bv_p(E))$\\ for $1\leq p < +\infty$} & $\Leftrightarrow$ & $f$ is constant\\[3mm] \midrule 
\parbox[c]{4.4cm}{\center $f \in H_{\text{bnd}}^\alpha(l^p(E),bv_q(E))$\\ for $\frac{p}{q}\leq \alpha \leq 1$} & $\Leftrightarrow$ & \parbox[c]{11cm}{\center $f$ is H\"older continuous on bounded sets with exponent $\alpha$} \\[3mm] \midrule

\parbox[c]{4.4cm}{\center $f \in H_{\text{bnd}}^\alpha(l^p(E),bv_q(E))$\\ for $0<\alpha \leq \frac{p}{q}\leq 1$} & $\Leftarrow$ & \parbox[c]{11cm}{\center $f$ is H\"older continuous on bounded sets with exponent $\alpha$ and there exits a ball $\ball_E(0,R)$ on which $f$ is H\"older continuous with exponent $\frac{p}{q}$} \\[3mm] \midrule

\parbox[c]{4.4cm}{\center $f \in H_{\text{bnd}}^\alpha(l^p(E),bv_q(E))$\\ for $0<\alpha\leq 1 < \frac{p}{q}$} & $\Leftrightarrow$ & \parbox[c]{11cm}{\center ??} \\[3mm] \midrule

\parbox[c]{4.4cm}{\center $f \in H_{\text{bnd}}^\alpha(bv_p(E),bv_q(E))$\\ for $1\leq q < p <+\infty$} & $\Leftrightarrow$ & $f$ is constant \\[3mm] \midrule 

\parbox[c]{4.4cm}{\center $f \in H_{\text{bnd}}^1(bv_1(\hat E),bv_p(\hat E))$\\ for $1\leq p < +\infty$} & $\Leftarrow$ & \parbox[c]{10.5cm}{\center $f$ is Fr\'echet differentiable with the derivative $f'$ being H\"older continuous on bounded sets with exponent $\frac{1}{p}$ } \\[4mm] \midrule 

\parbox[c]{4.4cm}{\center $f \in H_{\text{bnd}}^1(bv_1(\mathbb R),bv_1(\mathbb R))$} & $\Leftrightarrow$ & \parbox[c]{11cm}{\center the derivative of $f$ is Lipschitz continuous on bounded sets} \\[2mm] \midrule 

\parbox[c]{4.4cm}{\center $f \in H_{\text{bnd}}^\alpha (bv_p(E),bv_q(E))$\\ for $1< p \leq q<+\infty$} & $\Leftarrow$ & \parbox[c]{10.5cm}{\center the map $u \mapsto f(u)-f(0)$ is continuous and $\mathbb R$-linear} \\[3.5mm] \midrule 

\parbox[c]{4.4cm}{\center $f \in H_{\text{bnd}}^\alpha(bv_p(E),c_0(E))$\\ for $1 \leq p <+\infty$} & $\Leftrightarrow$ & $f$ is the zero map \\[3mm] \midrule 
\parbox[c]{4.4cm}{\center $f \in H_{\text{bnd}}^\alpha(bv_p(E),l^q(E))$\\ for $1 \leq p,q <+\infty$} & $\Leftrightarrow$ & $f$ is the zero map \\[3mm] \midrule 
$f \in H_{\text{bnd}}^\alpha(bv_1(\hat E),l^\infty(\hat E))$ & $\Leftrightarrow$ & $f$ is H\"older continuous on bounded sets with exponent $\alpha$\\ \midrule
\parbox[c]{4.4cm}{\center $f \in H_{\text{bnd}}^\alpha(bv_p(E),l^\infty(E))$\\ for $1<p<+\infty$} & $\Leftrightarrow$ & $f$ is constant \\[3mm] \midrule
$f \in H_{\text{bnd}}^\alpha(bv_1(\hat E),c(\hat E))$ & $\Leftrightarrow$ & $f$ is H\"older continuous on bounded sets with exponent $\alpha$ \\ \midrule
\parbox[c]{4.4cm}{\center $f \in H_{\text{bnd}}^\alpha(bv_p(E),c(E))$\\ for $1<p<+\infty$} & $\Leftrightarrow$ & $f$ is constant \\[2mm]
\bottomrule
\end{tabular}
\vspace{1.5mm}
\caption{H\"older continuity on bounded sets of the composition operator $C_f$.}
\end{center}
\end{table}
\end{small}

\clearpage

\begin{small}
\begin{table}[hbt!]
\begin{center}
\begin{tabular}{@{}ccc@{}}
\toprule
\parbox[c]{4.4cm}{\center $f \in H^\alpha(c_0(E),bv_p(E))$\\ for $1 \leq p < +\infty$} & $\Leftrightarrow$ & $f$ is constant\\[3mm] \midrule 
\parbox[c]{4.4cm}{\center $f \in H^\alpha(c(E),bv_p(E))$\\ for $1\leq p < +\infty$} & $\Leftrightarrow$ & $f$ is constant\\[3mm] \midrule 
\parbox[c]{4.4cm}{\center $f \in H^\alpha(l^\infty(E),bv_p(E))$\\ for $1\leq p < +\infty$} & $\Leftrightarrow$ & $f$ is constant\\[3mm] \midrule 
\parbox[c]{4.4cm}{\center $f \in H^\alpha(l^p(E),bv_q(E))$\\ for $\frac{p}{q}\leq \alpha \leq 1$} & $\Leftrightarrow$ & \parbox[c]{11cm}{\center $f$ is H\"older continuous with exponent $\alpha$} \\[3mm] \midrule

\parbox[c]{4.4cm}{\center $f \in H^\alpha(l^p(E),bv_q(E))$\\ for $1 \leq p,q < + \infty$\\ and $\alpha \in (0,1] \cap (0,\frac{p}{q})$} & $\Leftrightarrow$ & \parbox[c]{11cm}{\center $f$ is constant} \\[3mm] \midrule

\parbox[c]{4.4cm}{\center $f \in H^\alpha(bv_p(E),bv_q(E))$\\ for $1\leq q < p <+\infty$} & $\Leftrightarrow$ & $f$ is constant \\[3mm] \midrule 

\parbox[c]{4.4cm}{\center $f \in H^1(bv_1(\hat E),bv_p(\hat E))$\\ for $1\leq p < +\infty$} & $\Leftrightarrow$ & \parbox[c]{10.5cm}{\center the map $u \mapsto f(u)-f(0)$ is continuous and $\mathbb R$-linear} \\[3mm] \midrule 

\parbox[c]{4.4cm}{\center $f \in H^\alpha(bv_1(\hat E),bv_p(\hat E))$\\ for $1\leq p < +\infty$\\and $\alpha \in (0,1)$} & $\Leftrightarrow$ & \parbox[c]{10.5cm}{\center $f$ is constant} \\[3mm] \midrule 

\parbox[c]{4.4cm}{\center $f \in H^1(bv_p(E),bv_q(E))$\\ for $1< p \leq q<+\infty$} & $\Leftrightarrow$ & \parbox[c]{12cm}{\center the map $u \mapsto f(u)-f(0)$ is continuous and $\mathbb R$-linear} \\[3mm] \midrule 

\parbox[c]{4.4cm}{\center $f \in H^\alpha(bv_p(E),bv_q(E))$\\ for $1< p \leq q<+\infty$\\ and $\alpha \in (0,1)$} & $\Leftrightarrow$ & \parbox[c]{12cm}{\center $f$ is constant} \\[3mm] \midrule 
 
\parbox[c]{4.4cm}{\center $f \in H^\alpha(bv_p(E),c_0(E))$\\ for $1 \leq p <+\infty$} & $\Leftrightarrow$ & $f$ is the zero map \\[3mm] \midrule 
\parbox[c]{4.4cm}{\center $f \in H^\alpha(bv_p(E),l^q(E))$\\ for $1 \leq p,q <+\infty$} & $\Leftrightarrow$ & $f$ is the zero map \\[3mm] \midrule 
$f \in H^\alpha(bv_1(\hat E),l^\infty(\hat E))$ & $\Leftrightarrow$ & $f$ is  H\"older continuous with exponent $\alpha$\\ \midrule
\parbox[c]{4.4cm}{\center $f \in H^\alpha(bv_p(E),l^\infty(E))$\\ for $1<p<+\infty$} & $\Leftrightarrow$ & $f$ is constant \\[3mm] \midrule
$f \in H^\alpha(bv_1(\hat E),c(\hat E))$ & $\Leftrightarrow$ & $f$ is  H\"older continuous with exponent $\alpha$ \\ \midrule
\parbox[c]{4.4cm}{\center $f \in H^\alpha(bv_p(E),c(E))$\\ for $1<p<+\infty$} & $\Leftrightarrow$ & $f$ is constant \\[2mm]
\bottomrule
\end{tabular}
\vspace{1.5mm}
\caption{H\"older continuity of the composition operator $C_f$.}
\end{center}
\end{table}
\end{small}

\clearpage

\section{Compactness}

Finally, we characterize the compactness of composition operators. Since most of the work has been already done, we do not divide this part into subsections.

Before proceeding, let us recall that a (nonlinear) operator $F \colon X \to Y$ between two normed spaces $X$ and $Y$ is \emph{locally compact} if the image $F(B)$ of every bounded subset $B$ of $X$ is relatively compact in $Y$, or equivalently, is contained within a compact subset of $Y$. Similarly, $F$ is \emph{compact} if its entire range $F(X)$ is contained within a compact subset of $Y$. It is important to note that in this work we do not assume continuity when defining compactness or local compactness. A closely related concept is \emph{complete continuity}, which describes operators that are both pointwise continuous and locally compact.

It turns out that in our setting, all three properties -- compactness, local compactness, and complete continuity -- coincide. As we will show, locally compact composition operators in the sequence spaces we consider are generated solely by either a constant function or the zero map.

\begin{theorem}\label{thm:compact_c_l_bv_q}
Let $p \in [1,+\infty)$, and assume that $E$ is a Banach space when $p=1$, and a normed space when $p>1$. Additionally, let $Y \in \{c(E),l^\infty(E)\}$, and suppose that the composition operator $C_f$ maps $bv_p(E)$ into $Y$. Then, $C_f$ is locally compact if and only if $f$ is a constant map.
\end{theorem}

\begin{proof}
Clearly, it suffices to prove only the necessity part. Assume that $C_f \colon bv_p(E) \to Y$ is locally compact and suppose, for the sake of contradiction, that $f$ is not constant on $E$. Then, there exists a non-zero element $u \in E$ such that $f(u)\neq f(0)$. For each $n \in \mathbb N$ define the sequence $x_n:=(0,\ldots,0,u,0,0,\ldots)$, where the term $u$ is located at the $n$-th position. Notice that $\norm{x_n}_{bv_p}\leq 2\norm{u}$ for all $n \in \mathbb N$, so the set $B:=\dset{x_n}{n \in \mathbb N}$ is bounded in $bv_p(E)$. As the composition operator $C_f$ is locally compact, the image $C_f(B)$ is relatively compact in $Y$. In particular, the sequence $(C_f(x_n))_{n \in \mathbb N}$ has a subsequence that converges to some point $y \in Y$. It is easy to check that $y(m)=f(0)$ for all $m \in \mathbb N$. However, we also have $\norm{C_f(x_n)-y}_\infty = \norm{f(u)-f(0)}>0$ for all $n \in \mathbb N$. This implies that no subsequence of $(C_f(x_n))_{n \in \mathbb N}$ can converge to $y$, which is a contradiction. Therefore, $f$ must be constant on $E$.
\end{proof}

The situation remains unchanged if we consider composition operators acting between $l^p(E)$ and $bv_q(E)$.

\begin{theorem}
Let $p,q \in [1,+\infty)$. Moreover, assume that the composition operator $C_f$ maps $l^p(E)$ into $bv_q(E)$. Then, $C_f$ is locally compact if and only if $f$ is a constant map. 
\end{theorem}

\begin{proof}
The first part of the proof closely mirrors the argument used in the proof of Theorem~\ref{thm:compact_c_l_bv_q}. We include it here for the sake of completeness.

Assume that the composition operator $C_f \colon l^p(E) \to bv_q(E)$ is locally compact, and suppose that $f$ is not constant on $E$. Then, there exists a non-zero element $u \in E$ such that $f(u)\neq f(0)$. For each $n \in \mathbb N$ define the sequence $x_n:=(0,\ldots,0,u,0,0,\ldots)$, where $u$ appears in the $n$-th position. Note that $\norm{x_n}_{l^p}=\norm{u}$ for every $n \in \mathbb N$, so the set $B:=\dset{x_n}{n \in \mathbb N}$ is bounded in $l^p(E)$. Since $C_f$ is locally compact, the image $C_f(B)$ is relatively compact in $bv_q(E)$. In particular, the sequence $(C_f(x_n))_{n \in \mathbb N}$ has a subsequence $(C_f(x_{n_k}))_{k \in \mathbb N}$ that converges to some point $y \in bv_q(E)$. We will show that $y(m)=f(0)$ for all $m \in \mathbb N$. 

First, observe that for every $k \geq 2$ we have $\norm{f(0)-y(1)}\leq \norm{C_f(x_{n_k})-y}_{bv_q}$. Since the right-hand side tends to zero as $k \to +\infty$, it follows that $y(1)=f(0)$. Now, suppose that $y(m)=f(0)$ for all $m=1,\ldots,l$. Then, for $k \geq l+2$ we have 
\begin{align*}
& \norm{C_f(x_{n_k})-y}_{bv_q}^q\\
&\qquad \geq \sum_{m=1}^l\norm[\Big]{ \bigl[f(x_{n_k}(m+1))-y(m+1)\bigr] - \bigl[f(x_{n_k}(m))-y(m)\bigr]}^q\\
&\qquad=\sum_{m=1}^l\norm[\big]{y(m+1) - y(m)}^q = \norm{y(l+1)-f(0)}^q.
\end{align*}
Since the left-hand side tends to zero as $k \to +\infty$, we deduce that $y(l+1)=f(0)$. Therefore, by induction, $y(m)=f(0)$ for every $m \in \mathbb N$. However, for $k \geq 2$ we also have $\norm{C_f(x_{n_k})-y}_{bv_q} \geq 2^{\frac{1}{q}}\norm{f(u)-f(0)}>0$, which contradicts the fact that $(C_f(x_{n_k}))_{k \in \mathbb N}$ converges to $y$. This contradiction completes the argument. 
 \end{proof}

Combining the arguments above with the fact that the composition operator $C_f$ maps $bv_p(E)$ into $bv_q(E)$ for $1 \leq q < p < +\infty$ if and only if $f$ is a constant map (see \cite{BK}*{Corollary 19}), we obtain the following result.

\begin{theorem}\label{thm:compactnes_bv_p_bv_q}
Let $p,q \in [1,+\infty)$, and assume that $E$ is a Banach space when $p=1$, and a normed space when $p>1$. Assume further that the composition operator $C_f$ maps $bv_p(E)$ into $bv_q(E)$. Then, $C_f$ is locally compact if and only if $f$ is a constant map.
\end{theorem}

\subsection{Compactness: summary}

Let $E$ be a normed space over the real or complex field. If $E$ is a Banach space, we denote it by $\hat E$ to emphasise its completeness. For two sequence spaces $X$ and $Y$ we denote by $K(X,Y)$ the class of all maps $f\colon E \to E$ such that the generated composition operator $C_f \colon X \to Y$ is either locally compact, compact, or completely continuous. Then, taking into account the results established in the previous parts of the paper, along with the results on acting conditions and boundedness presented in~\cite{BK}*{Sections~3 and~4}, we summarize the discussion of this section in the following table.  

\clearpage

\begin{small}
\begin{table}[hbt!]
\begin{center}
\begin{tabular}{@{}ccc@{}}
\toprule
\parbox[c]{4cm}{\center $f \in K(c_0(E),bv_p(E))$\\ for $1 \leq p < +\infty$} & $\Leftrightarrow$ & \parbox[c]{4.6cm}{\center $f$ is constant}\\[3mm] \midrule 
\parbox[c]{4cm}{\center $f \in K(c(E),bv_p(E))$\\ for $1\leq p < +\infty$} & $\Leftrightarrow$ & $f$ is constant\\[3mm] \midrule 
\parbox[c]{4cm}{\center $f \in K(l^\infty(E),bv_p(E))$\\ for $1\leq p < +\infty$} & $\Leftrightarrow$ & $f$ is constant\\[3mm] \midrule 
\parbox[c]{4cm}{\center $f \in K(l^p(E),bv_q(E))$\\ for $1\leq p , q < +\infty$} & $\Leftrightarrow$ & $f$ is constant \\[3mm] \midrule
\parbox[c]{4cm}{\center $f \in K(bv_p(E),bv_q(E))$\\ for $1\leq q < p <+\infty$} & $\Leftrightarrow$ & $f$ is constant \\[3mm] \midrule 
\parbox[c]{4cm}{\center $f \in K(bv_1(\hat E),bv_p(\hat E))$\\ for $1\leq p < +\infty$} & $\Leftrightarrow$ & $f$ is constant \\[3mm] \midrule 
\parbox[c]{4cm}{\center $f \in K(bv_p(E),bv_q(E))$\\ for $1< p \leq q<+\infty$} & $\Leftrightarrow$ & $f$ is constant \\[3.5mm] \midrule 
\parbox[c]{4cm}{\center $f \in K(bv_p(E),c_0(E))$\\ for $1 \leq p <+\infty$} & $\Leftrightarrow$ & $f$ is the zero map \\[3mm] \midrule 
\parbox[c]{4cm}{\center $f \in K(bv_p(E),l^q(E))$\\ for $1 \leq p,q <+\infty$} & $\Leftrightarrow$ & $f$ is the zero map \\[3mm] \midrule 
$f \in K(bv_1(\hat E),l^\infty(\hat E))$ & $\Leftrightarrow$ & $f$ is constant\\ \midrule
\parbox[c]{4cm}{\center $f \in K(bv_p(E),l^\infty(E))$\\ for $1<p<+\infty$} & $\Leftrightarrow$ & $f$ is constant \\[3mm] \midrule
$f \in K(bv_1(\hat E),c(\hat E))$ & $\Leftrightarrow$ & $f$ is constant \\ \midrule
\parbox[c]{4.4cm}{\center $f \in K(bv_p(E),c(E))$\\ for $1<p<+\infty$} & $\Leftrightarrow$ & $f$ is constant \\[2mm]
\bottomrule
\end{tabular}
\vspace{1.5mm}
\caption{Compactness properties of the composition operator $C_f$.}
\end{center}
\end{table}
\end{small}

\clearpage

\section{Discussion and conclusions}

In this final section, following the approach in~\cite{BK}*{Section~5}, we compare our results with selected theorems from the parallel theory of composition operators in the $BV_p$-spaces of functions of bounded Wiener variation. Our discussion will be limited to operators acting within the same space and to the case $E:=\mathbb R$, as little is currently known about the setting in which $E$ is a general normed space and the domain and target $BV_p$-spaces correspond to different parameters~$p$.     

We begin by recalling the notion of Wiener variation. Let $p \in [1,+\infty)$ be fixed. The (possibly infinite) quantity
\[
 \var_p x:=\sup \sum_{i=1}^n\abs[\big]{x(t_{i})-x(t_{i-1})}^p,
\]
where the supremum is taken over all finite partitions $a=t_0<t_1<\cdots<t_n=b$ of $[a,b]$, is called the \emph{$p$-variation} (or, \emph{Wiener variation}) of the real-valued map $x$ over the interval $[a,b]$. A mapping $x \colon [a,b] \to \mathbb R$ is said to be of \emph{bounded $p$-variation} if $\var_p x<+\infty$. The space of all such maps, denoted by $BV_p[a,b]$, forms a normed space when endowed with the norm $\norm{x}_{BV_p}:=\abs{x(a)} + (\var_p x)^{1/p}$. For further information on Wiener variation, we refer the reader to~\cite{ABM} and~\cite{DN}.

In the functional setting, the composition operator $C_f$, generated by $f \colon \mathbb R \to \mathbb R$, is defined analogously to our case, with obvious adjustments. Specifically, to a given function $x \colon [a,b] \to \mathbb R$, the composition operator assigns the real-valued map $C_f(x)$ defined for $t \in [a,b]$ by $C_f(x)(t):=f(x(t))$. We will not introduce a separate symbol to distinguish composition operators in $bv_p$- and $BV_p$-spaces, as the context will always make the intended meaning clear.

At the heart of the theory of composition operators on spaces of functions of bounded variation lies the Josephy theorem, which states that $C_f$ maps $BV_p[a,b]$ into itself if and only if the function $f \colon \mathbb R \to\mathbb R$ is Lipschitz continuous on bounded subsets of $\mathbb R$ (see, for example, \cite{ABM}*{Theorem~5.12 and the subsequent discussion},~\cite{DN}*{Corollary~6.36} and~\cite{reinwand}*{Proposition~5.1.1}).

One of the most long-standing and important open questions in this context was whether the conditions ensuring that composition operators act on $BV_p$-spaces also imply their pointwise continuity. As noted in~\cite{BGK}, a positive answer for the case $p=1$ was, in fact, established by Morse as early as 1937. His proof, however, was highly technical and rather complex. Simpler proofs were later provided by Ma\'ckowiak in~\cite{M} and Reinwand in~\cite{SR} (see also~\cite{reinwand}*{Theorem~5.1.24}). 

Our findings on pointwise continuity of composition operators in $bv_p$-spaces are therefore fully consistent with the classical framework of functions of bounded variation. Moreover, they remain valid in the more general setting of composition operators acting between $bv_p(\mathbb{R})$ and $bv_q(\mathbb{R})$ for any $p,q \in [1, +\infty)$ -- see Theorem~\ref{thm:cont_bv_p_bv_q}. For completeness, we note that the question of whether every composition operator  $C_f \colon BV_p[a,b] \to BV_p[a,b]$ is continuous remains open.

The results concerning locally uniformly continuous and uniformly continuous composition operators are also consistent across both settings. In particular, a composition operator $C_f \colon bv_p(\mathbb{R}) \to bv_p(\mathbb{R})$ or $C_f \colon BV_p[a,b] \to BV_p[a,b]$ is locally uniformly continuous for $p=1$ if and only if its generator $f\colon \mathbb R \to\mathbb R$ is continuously differentiable (see Theorem~\ref{thm:LUC_C1} and~\cite{reinwand}*{Theorem~5.1.23}). Furthermore, for any $p \in [1,+\infty)$, these operators are uniformly continuous if and only if their generators are affine functions (see Corollary~\ref{cor:UC_p_q_R} and~\cite{reinwand}*{Theorem~5.1.22}).

The equivalence between the continuous differentiability of a function $f$ and the local uniform continuity of the associated composition operator $C_f$ also holds for operators acting on the space $BV_p[a,b]$ for $p>1$ (see \cite{reinwand}*{Theorem~5.1.22}). In contrast, whether a similar characterization applies to $bv_p$-spaces for $p>1$ remains an open question. The main difficulty in addressing this problem lies in the fact that, unlike functions of bounded $p$-variation which are always bounded, the space $bv_p(\mathbb{R})$ contains unbounded sequences when $p>1$. 

It is also worth noting that, within the framework of $BV_p$-spaces, no analogue of Theorem~\ref{thm:LUC_bv_1_bv_p_R} has yet been established. Nevertheless, we believe that the method developed in the proof of Theorem~\ref{thm:LUC_bv_1_bv_p_R} could be adapted to provide sufficient conditions for the local uniform continuity of composition operators between different spaces of functions of bounded Wiener variation.

Thus far, the study of composition operators on spaces of functions of bounded Wiener variation has focused primarily on their Lipschitz continuity, with comparatively little attention given to the H\"older condition. For this reason, the discussion here will be restricted to the Lipschitz case. The analysis of the (global) Lipschitz continuity for composition operators in both $bv_p$- and $BV_p$-spaces is relatively straightforward. Since these operators are obviously uniformly continuous, the associated degeneracy phenomena that arise are consistent across both the functional and sequential frameworks (see Corollary~\ref{cor:UC_p_q_R} and~\cite{ABM}*{Theorem~5.47} or~\cite{reinwand}*{Theorem~5.1.22}). 

A similar agreement holds for Lipschitz continuity on bounded sets when $p = 1$. In particular, in both settings, a composition operator $C_f$ is Lipschitz continuous on bounded subsets if and only if the derivative $f'$ of its generator $f \colon \mathbb R \to \mathbb R$ is Lipschitz continuous on bounded sets (see Theorem~\ref{thm:LH1_bv_1_bv_p_ver2} and~\cite{reinwand}*{Theorem~5.1.21}). This characterization also extends to operators acting on $BV_p[a,b]$ when $p > 1$ (see~\cite{reinwand}*{Theorem~5.1.21}). However, whether the same equivalence holds for composition operators on $bv_p(\mathbb{R})$ remains an open question. As in the case of local uniform continuity, the difficulty here lies in the fact that the space $bv_p(\mathbb{R})$ contains unbounded sequences for $p > 1$.

Finally, it is worth noting that the problem of (local) compactness for composition operators in both $bv_p(\mathbb{R})$ and $BV_p[a,b]$ is fully understood: in either setting, such operators must necessarily be constant, regardless of the choice of the parameter $p \in [1,+\infty)$ (see Theorem~\ref{thm:compactnes_bv_p_bv_q} and~\cite{reinwand}*{Theorem~5.1.20}).

\vskip .5cm

\textit{Declarations:
\begin{itemize}
\item On behalf of all authors, the corresponding author states that there is no conflict of interest.  
\item No funding was received to assist with the preparation of this manuscript.
\item The authors have no relevant financial or non-financial interests to disclose.
\end{itemize}}


\begin{bibdiv}
\begin{biblist}

\bib{ABM}{book}{
  title={Bounded variation and around},
  author={Appell, J.},
  author={Bana\'s, J.},
  author={Merentes, N.},
  date={2014},
  series={De Gruyter Studies in Nonlinear Analysis and Applications, vol. 17},
  publisher={De Gruyter},
  address={Berlin},
}

\bib{ABKR}{book}{
  title={BV type spaces with applications},
  author={Appell, J.},
  author={Bugajewska, D.},
  author={Kasprzak, P.},
	author={Reinwand, S.},
   series={Lecture Notes in Nonlinear Analysis},
   volume={20},
   publisher={Juliusz Schauder Center for Nonlinear Studies, Toru\'{n}},
   date={2021},
}

\bib{AZ}{book}{
   author={Appell, J.},
   author={Zabrejko, . P.},
   title={Nonlinear superposition operators},
   series={Cambridge Tracts in Mathematics},
   volume={95},
   publisher={Cambridge University Press, Cambridge},
   date={1990},
}

\bib{BA}{article}{
   author={Ba\c{s}ar, F.},
   author={Altay, B.},
   title={On the space of sequences of $p$-bounded variation and related
   matrix mappings},
   language={English, with English and Ukrainian summaries},
   journal={Ukra\"{\i}n. Mat. Zh.},
   volume={55},
   date={2003},
   number={1},
   pages={108--118},
   translation={
      journal={Ukrainian Math. J.},
      volume={55},
      date={2003},
      number={1},
      pages={136--147},
   },
}

\bib{BAM}{article}{
    AUTHOR = {Ba\c{s}ar, F.},
		AUTHOR = {Altay, B.},
		AUTHOR = {Mursaleen, M.},
     TITLE = {Some generalizations of the space {$bv_p$} of {$p$}-bounded
              variation sequences},
   JOURNAL = {Nonlinear Anal., Theory Methods Appl., Ser. A, Theory Methods},
    VOLUME = {68},
      YEAR = {2008},
    NUMBER = {2},
     PAGES = {273--287},
}


\bib{BBK}{book}{
   author={Borkowski, M.},
   author={Bugajewska, D.},
   author={Kasprzak, P.},
   title={Selected topics in nonlinear analysis},
   series={Lecture Notes in Nonlinear Analysis},
   volume={19},
   publisher={Juliusz Schauder Center for Nonlinear Studies, Toru\'{n}},
   date={2021},
}

\bib{BK}{article}{
   author={Bugajewska, D.},
   author={Kasprzak, P.},
   title={Nonlinear composition operators in $bv_p$-spaces\textup: acting conditions and boundedness},
	 pages={submitted},
	}

\bib{BGK}{article}{
   author={Bugajewski, D.}, 
	author={Gulgowski, J.},
	author={Kasprzak, P.},
	title={ On continuity and compactness of some nonlinear operators in the spaces of functions of bounded variation},
  journal={Ann. Mat. Pura Appl.},
	date={2016},
	volume={195},
	pages={1513--1530},
	}


\bib{DZ}{article}{
   author={Dedagich, F.},
   author={Zabre\u{\i}ko, P. P.},
   title={On superposition operators in $l_p$ spaces},
   language={Russian},
   journal={Sibirsk. Mat. Zh.},
   volume={28},
   date={1987},
   number={1},
   pages={86--98},
}

\bib{DL}{book}{
  author = {R. A. DeVore},
	author = {G. G. Lorentz},
  title = {Constructive Approximation},
  publisher = {Springer-Verlag},
  date = {1993},
  series = {Grundlehren der mathematischen Wissenschaften},
  volume = {303},
}

\bib{DN}{book}{
   author={Dudley, R. M.},
   author={Norvai\v{s}a, R.},
   title={Concrete functional calculus},
   series={Springer Monographs in Mathematics},
   publisher={Springer, New York},
   date={2011},
}

\bib{HLP}{book}{
   author={Hardy, G. H.},
   author={Littlewood, J. E.},
   author={P\'{o}lya, G.},
   title={Inequalities},
   publisher={Cambridge, at the University Press},
   date={1934},
}



\bib{KFA}{article}{
   author={Karami, S.},
   author={Fathi, J.},
	 author={Ahmadi, A.},
	 date={2023},
	 volume={14},
	 number={7},
   title={The results of the superposition operator on sequence space $bv_p$},
   journal={Int. J. Nonlinear Anal. Appl. },
   pages={309--320},
}

\bib{K}{article}{
    AUTHOR = {Kiri\c{s}ci, M.},
     TITLE = {The sequence space $bv$ and some applications},
   JOURNAL = {Math. \AE terna},
    VOLUME = {4},
      YEAR = {2014},
    NUMBER = {3--4},
     PAGES = {207--223},
}


\bib{LF}{article}{
    AUTHOR = {Lashkaripour, R.},
		AUTHOR = {Fathi, J.},
     TITLE = {Norms of matrix operators on {$bv_p$}},
   JOURNAL = {J. Math. Inequal.},
    VOLUME = {6},
      YEAR = {2012},
    NUMBER = {4},
     PAGES = {589--592},
}

\bib{M}{article}{
  author={Ma\'ckowiak, P.},
  title={On the continuity of superposition operators in the 
space of functions of bounded variation},
	journal={Aequationes Math.},
	volume={91},
	number={4},
	date={2017},
  pages={759--777},
}

\bib{MatMis}{article}{
  author={Matkowski, J.},
  author={Mi\'{s}, J.},
  title={On a characterization of Lipschitzian operators of substitution in the space ${\rm BV}\langle a$, $b\rangle $},
  journal={Math. Nachr.},
  volume={117},
  date={1984},
  pages={155--159},
}

\bib{Mo}{article}{
  author={Morse, A. P.},
  title={Convergence in variation and related topics},
	journal={Trans. Amer. Math. Soc.},
	volume={41},
	number={1},
	date={1937},
  pages={48--83},
	}


\bib{SR}{article}{
  author={Reinwand, S.},
  title={Types of convergence which preserve continuity},
	journal={Real Anal. Exchange},
	volume={45},
	number={1},
	date={2020},
  pages={173--204},
	}


\bib{reinwand}{book}{
  title={Functions of bounded variation\textup: Theory, methods, applications},
  author={Reinwand, S.},
  year={2021},
  publisher={Cuvillier Verlag}
}

\bib{small}{book}{
   author={Small, Ch. G.},
   title={Functional equations and how to solve them},
   series={Problem Books in Mathematics},
   publisher={Springer, New York},
   date={2007},
}



\end{biblist}
\end{bibdiv}
\end{document}